\documentclass{article}

\usepackage[british]{babel}

\usepackage[T1]{fontenc}
\usepackage{amsbsy}
\usepackage{amsthm}

\theoremstyle{plain}

\usepackage{latexsym}
\usepackage{amsmath}
\usepackage{amsfonts}
\usepackage{amssymb}
\usepackage{psfrag}
\usepackage{graphicx}

\newcommand{\be}{\begin{equation}}
\newcommand{\ee}{\end{equation}}

\newcommand{\bem}{\begin{multline}}
\newcommand{\eem}{\end{multline}}

\newcommand{\bml}{\begin{multline*}}
\newcommand{\eml}{\end{multline*}}

\newcommand{\beg}{\begin{gather}}
\newcommand{\eeg}{\end{gather}}

\theoremstyle{definition}
\newtheorem{definition}{Definition}

\newtheorem{exm}{Example}
\newtheorem{lemma}{Lemma}
\newtheorem{Corollary}{Corollary}

\newtheorem{thm}{Theorem}

\begin{document}

\title{On interpolation problem for multidimensional harmonizable stable sequences with noise observations}
\author{Mikhail Moklyachuk$^1$
\footnote{*Corresponding email:Moklyachuk@gmail.com%
 \newline 1 Department of Probability Theory, Statistics and Actuarial
Mathematics, Taras Shevchenko National University of Kyiv, Kyiv, Ukraine%
%
 %\newline 1E-mail: youremail@gmail.com%
 %\newline 2E-mail: youremail@yahoo.com
 }
\,
 }
\maketitle

%\artinfo{XX XXXXXX XXXX}{XX XXXXXX XXXX}

\abstract{ We consider the problem of optimal linear estimation of the functional
$$A_N \vec{\xi} =\sum_{j = 0}^{N} (\vec{a}(j))^{\top} \vec{\xi}(j)$$
that depends on the
unknown values  $\vec{\xi}(j),j=0,1,\dots,N,$  of a vector-valued harmonizable symmetric $\alpha$-stable random sequence
$\vec{\xi}(j)=\left \{ \xi_ {k} (j) \right \}_{k = 1} ^ {T}$,
from observations of the sequence $\vec{\xi}(j)+\vec{\eta}(j)$ at points  $j\in\mathbb Z\setminus\{0,1,\dots,N\}$.
We consider the problem for mutually independent
vector-valued harmonizable symmetric $\alpha$-stable random sequences $\vec{\xi}(j)=\left \{ \xi_ {k} (j) \right \}_{k = 1} ^ {T}$
and
 $\vec{\eta}(j)=\left \{ \xi_ {k} (j) \right \}_{k = 1} ^ {T}$
 which have
absolutely continuous spectral measures and the spectral densities $f(\theta)$ and $g(\theta)$
satisfying the minimality condition.
}

\vspace{2ex}

\vspace{2ex}

\section{Introduction}

The problem of estimation of the unknown values of harmonizable random sequences and processes were investigated in papers by
Cambanis (1983),
 Cambanis and Soltani (1984), Hosoya (1982).
The interpolation problem for harmonizable symmetric $\alpha$-stable random sequences were investigated in papers by
Weron (1985) and Pourahmadi (1984).

Basic results concerning estimation of the unknown (missed) values of stochastic processes are based on the assumption that spectral densities of processes are exactly known. In practice, however, complete information on the spectral densities is impossible in most cases.
In such situations one finds parametric or nonparametric estimates of the unknown spectral densities.
Then the classical estimation method is applied under the assumption that the estimated densities are true.
This procedure can result in significant increasing of the value of error as Vastola and Poor (1983)  have demonstrated with the help of some examples.
 This is a reason to search estimates which are optimal for all densities from a class of admissible spectral densities.
 These estimates are called minimax since they minimize the maximal value of the error.
 A survey of results (till 1985) in minimax (robust) methods of data processing  can be found in the paper by Kassam and Poor (1985).
 The paper by Grenander (1957) should be marked as the first one where the minimax extrapolation problem for stationary processes was formulated and solved.
Later Franke and Poor (Franke and Poor, 1984; Franke, 1985) applied the convex optimization methods for investigation the minimax-robust extrapolation and interpolation problems.
In papers by Moklyachuk (1994 -- 2008) of the minimax-robust extrapolation, interpolation and filtering  problems are studied for stationary processes.
The book by Moklyachuk and Masyutka (2012) are dedicated to minimax-robust extrapolation, interpolation and filtering  problems for vector-valued stationary processes and sequences.
  In the book by  Moklyachuk and Golichenko (2016) the minimax-robust extrapolation, interpolation and filtering  problems for periodically correlated processes are investigated.
 The minimax-robust extrapolation, interpolation and filtering  problems for stochastic processes with $n$th stationary increments are investigated by
 Luz and Moklyachuk (2019 -- 2024).
 In papers by Moklyachuk and Ostapenko (2015, 2016)  minimax-robust interpolation  problems are studied for harmonizable stable sequences.

In this paper the problem of optimal estimation is investigated for the linear functional
$$A_N \vec{\xi} =\sum_{j = 0}^{N} (\vec{a}(j))^{\top} \vec{\xi}(j)
$$
that depends on the
unknown values  $\vec{\xi}(j),j=0,1,\dots,N,$  of a vector-valued harmonizable symmetric $\alpha$-stable random sequence $\vec{\xi}(j)=\left \{ \xi_ {k} (j) \right \}_{k = 1} ^ {T}$,
from observations of the sequence $\vec{\xi}(j)+\vec{\eta}(j)$ at points  $j\in\mathbb Z\setminus\{0,1,\dots,N\}$
where $\vec{\xi}(j)$ and $\vec{\eta}(j)$ are mutually independent harmonizable symmetric $\alpha$-stable random sequences
which have the spectral densities $f(\theta)$ and $g(\theta)$ satisfying the minimality condition.

The problem is investigated under the condition of spectral certainty as well as under the condition of spectral uncertainty.
 Formulas for calculation the value of the error and spectral characteristic of the optimal linear
estimate of the functional are derived under the condition of spectral certainty where spectral density of the sequence is exactly known.
 In the case where spectral density of
the sequence is not exactly known  sets of admissible spectral densities is available, relations which determine least favorable densities and the minimax-robust spectral characteristics for different classes of spectral densities are found.

 \section{Harmonizable symmetric $\alpha$-stable random sequence}

 \begin{definition}(symmetric $\alpha$-stable random variable)
 A real-valued random variable $\xi$ is said to be symmetric $\alpha$-stable, $S\alpha S$, if its characteristic function has the form $E exp(it\xi) = exp(-c|t|^{\alpha})$ for some $c \geq 0$ and $0 < \alpha \leq 2.$
 The real random variables $\xi_1,\xi_2,\dots,\xi_n$ are jointly $S\alpha S$ if all linear combinations
 $\sum_{k=1}^{n}a_k\xi_k$ are $S\alpha S$, or, equivalently, if the characteristic function of the random vector  $\vec{\xi}=(\xi_1,\dots,\xi_n)$ is of the form
 $$\phi_{\vec{\xi}}(\vec{t}) = E exp\left(i \sum_{k=1}^{n} t_k \xi_k\right) = exp\left\{-\int_{S_n}\left|\sum_{k=1}^{n}t_k x_k\right|^{\alpha}d \Gamma_{\vec{\xi}}(\vec{x})\right\},$$
  where $\vec{t}=(t_1,\dots, t_n),  t_1,\dots, t_n$ are real numbers and $\Gamma_{\vec{\xi}}(\vec{x})$ is a symmetric measure defined on the unit sphere $S_n \in R^n$, called the spectral measure of the random vector  $\vec{\xi}=(\xi_1,\dots,\xi_n)$. There is a one-to-one correspondence between the distribution of  $\vec{\xi}$ and its spectral measure $\Gamma_{\vec{\xi}}(\vec{x})$  (Cambanis, 1983).
 \end{definition}

   For real-valued  jointly $S\alpha S$ random variables $\xi,\eta$ with $1<\alpha\leq 2$ the covariation of $\xi,\eta$  is defined by
  $$[\xi,\eta]_{\alpha} = \int_{S_2} (x)(y)^{<\alpha-1>} d \Gamma_{\xi,\eta}(x, y),$$

 \noindent where $(y)^{<\beta>} = |y|^{\beta - 1} {y} $.

 For jointly $S{\alpha}S$ random variables $\xi= \xi_1 + i\xi_2$ and $\eta= \eta_1 + i \eta_2$ the covariation of $\xi$ with $\eta$ is defined as (Cambanis, 1983)
 $$[\xi,\eta]_{\alpha} = \int_{S_4} (x_1 + i x_2)(y_1 + i y_2)^{<\alpha-1>} d \Gamma_{\xi_1,\xi_2,\eta_1,\eta_2}(x_1, x_2, y_1, y_2),$$

 \noindent where $z^{<\beta>} = |z|^{\beta - 1} \bar{z} $ for a complex number $z$ and $\beta > 0.$

 The covariation in general is not symmetric and linear on second argument and for $\xi, \xi_1, \xi_2, \eta$ jointly $S\alpha S$  has the following properties (Cambanis, 1983, Weron, 1985).

 %\begin{itemize}

  \begin{equation}\label{eq:stable1}
 [\xi_1 + \xi_2, \eta]_{\alpha} = [\xi_1, \eta]_{\alpha} + [\xi_2, \eta]_{\alpha},
 \end{equation}

  \begin{equation}\label{eq:stable2}
 [a\,\xi, b\,\eta]_{\alpha} = a(b)^{\alpha-1}\, [\xi, \eta]_{\alpha},
 \end{equation}

   \begin{equation}\label{eq:stable3}
 [\xi, \eta]_{\alpha} = 0 \quad \text{if $\xi$ and $\eta$ are independent},
 \end{equation}

   \begin{equation}\label{eq:stable4}
 [\xi, \eta_1+\eta_2]_{\alpha} = [\xi, \eta_1]_{\alpha} + [\xi, \eta_2]_{\alpha}\quad\text{if $\eta_1$ and $\eta_2$ are independent},
 \end{equation}

 \begin{equation}\label{eq:inequality5}
 |[\xi, \eta]_{\alpha}| \leq ||\xi||_{\alpha} ||\eta||_{\alpha}^{\alpha - 1},
 \end{equation}
the functional
  \begin{equation}\label{eq:stable6}
  ||\xi||_{\alpha} = [\xi, \xi]_{\alpha}^{1/\alpha}
   \end{equation}

 \noindent  is a norm in a linear space of $S\alpha S$ random variables which is equivalent to convergence in probability.

  \begin{equation}\label{eq:stable6}
  \|\sum_{k=1}^n\xi_k\|_{\alpha}^{\alpha} =\sum_{k=1}^n \|\xi_k\|_{\alpha}^{\alpha}
  \end{equation}
 \noindent when $\xi_1,\dots,\xi_n$ are independent,

 \noindent the mapping
    \begin{equation}\label{eq:stable7}
  \xi\to [\xi,\eta]_{\alpha}
  \end{equation}

 \noindent is a bounded linear functional with the norm $\|\eta\|_{\alpha}^{\alpha-1}$ on the linear space of $S\alpha S$ random variables, and every bounded linear functional on such a space
 is of this form for some $\eta$.

It should be noted that $||\cdot||_\alpha$ is not necessarily the usual $L^{\alpha}$ norm.

 Here is the simplest properties of the function $z^{<\beta>}.$
\begin{lemma}
 Let $z, x, y $ be complex numbers, $\beta > 0$. Then the following properties holds true:
 \begin{itemize}
 \item $|z|^{<\beta>} = z \cdot z^{<\beta - 1>},$
 \item $\left||z|^{<\beta>}\right| = \left|z\right|^{<\beta>},$
 \item if $z^{<\beta>} = v$, thet $z = v^{<1/\beta>} = |v|^{(1-\beta)/\beta}\bar{v},$
 \item $z^{<1>} = \bar{z},$
 \item if $z \neq 0$, then $z^{<\alpha>} z^{<\beta>} = \frac{\bar{z}}{|z|} z^{<\alpha + \beta>},$
 \item if $z \neq 0$, then $\frac{z^{<\alpha>}}{z^{<\beta>}} = \frac{z}{|z|} z^{<\alpha - \beta>},$
 \item $(c z)^{<\alpha>} = c^{\alpha} z^{<\alpha>}, c \in \mathbb{R},$
 \item $(z^{<\alpha>})^{<\beta>} = {\bar{z}}^{<\alpha \beta>},$
 \item $(xy)^{<\alpha>} = x^{<\alpha>} y^{<\alpha>},$
 \item $(z^{\alpha})^{<\beta>} = (z^{<\beta>})^{\alpha},$
 \item $(z^{<\alpha>})^{\beta} = (z^{\beta})^{<\alpha>},$
 \item $|z^{<\alpha>}|^{\beta} = |z|^{\alpha \beta},$
 \item $(x + y)^{<\alpha>} = \bar{x}|x + y|^{\alpha - 1} + \bar{y}|x + y|^{\alpha - 1}.$
\end{itemize}
\end{lemma}

 \begin{definition}[symmetric $\alpha$-stable stochastic sequence]
A stochastic sequence $\{\xi(n),n\in\mathbb Z\}$ is called symmetric $\alpha$-stable, $S\alpha S$,
if all  linear combinations  $\sum_{m=1}^{l}a_{m}\xi(n_m)$ are $S\alpha S$ random variables.

A vector-valued stochastic sequence $\vec{\xi}(n)=\left \{ \xi_ {k} (n) \right \}_{k = 1} ^ {T},\,n\in \mathbb{Z}$,
is called symmetric $\alpha$-stable, $S\alpha S$, stochastic sequence,
if all  linear combinations \linebreak
 $\sum_{k=1}^{T}\sum_{m=1}^{l}a_{m_k}\xi_k(n_m)$ are $S\alpha S$ random variables.
\end{definition}

Let $Z =\{Z(t): -\infty < t < \infty\}$ be a complex $S{\alpha}S$ process with independent increments. The spectral measure of the process $Z$ is defined as $\mu\{(s, t]\} =\|Z(t) - Z(s)\|_{\alpha}^{\alpha}.$

The integrals $\int a(t)dZ(t)$ can be defined for all $a(t) \in L^{\alpha}(\mu)$ with properties for all $a \in L^{\alpha}(\mu),b \in L^{\alpha}(\mu)$ (see Cambanis, 1983; Cambanis and Soltani, 1984; Hosoya, 1982):
\begin{equation}\label{eq:norm_equality11}
 \left\|\int a(t) d Z(t)\right\|^{\alpha}_{\alpha} = \int |a(t)|^{\alpha}d \mu(t),
 \end{equation}
  \begin{equation}\label{eq:norm_equality12}
 \left[\int a(t) d Z(t), \int b(t) d Z(t)\right]_\alpha = \int a(t) (b(t))^{<\alpha - 1>} d \mu(t),
\end{equation}
and for vector-valued functions
$\vec{a}(t)=\left \{ a_k(t) \right \}_{k = 1} ^ {T},\,a_k(t)\in L^{\alpha}(\mu)$,
$\vec{b}(t)=\left \{ b_k(t) \right \}_{k = 1} ^ {T},\,b_k(t)\in L^{\alpha}(\mu)$,
and
$\vec{Z}(t)=\left \{ Z_k(t) \right \}_{k = 1} ^ {T}$,
  \begin{equation}\label{eq:norm_equality13}
 \left[\int (\vec{a}(t))^{\top}  d \vec{Z}(t), \int (\vec{b}(t))^{\top}  d \vec{Z}(t)\right]_\alpha = \int (\vec{a}(t))^{\top} d \mu(t) (\vec{b}(t))^{<\alpha - 1>},
\end{equation}
where $\mu$ is the matrix-valued spectral measure corresponding to the process $\vec{Z}(t)$.

\begin{definition}[harmonizable symmetric $\alpha$-stable sequence]
 A symmetric $\alpha$-stable, $S\alpha S$, stochastic sequence $\{\xi(n),n\in\mathbb Z\}$ is said to be harmonizable, $HS{\alpha}S$, if there exists a $S\alpha S$ process
 $Z = \{Z(\theta): \theta \in [-\pi, \pi]\}$ with independent increments and a finite spectral measure $\mu$ such that sequence $\xi(n)$ has the spectral representation
 $$\xi(n) = \int_{-\pi}^{\pi}e^{in\theta}dZ(\theta), \quad n \in \mathbb{Z},$$
 \noindent and the covariation has the representation
 $$[\xi(n), \xi(m)]_\alpha = \int_{-\pi}^{\pi}e^{i(n-m)\theta}d \mu(\theta), \quad m, n \in \mathbb{Z}.$$

 A vector-valued symmetric $\alpha$-stable, $S\alpha S$, stochastic sequence
 $\vec{\xi}(n)=\left \{ \xi_k(n) \right \}_{k = 0} ^ {T-1},\,n\in \mathbb{Z}$,
 is said to be harmonizable, $HS{\alpha}S$, if there exists a vector-valued  $S\alpha S$ process
 $\vec{Z}(\theta) = \left\{{Z}_k(\theta)\right\}_{k = 0} ^ {T-1}, \theta \in [-\pi, \pi)$
 with independent increments and a finite matrix-valued spectral measure
 $\mu$ such that sequence
 $\vec{\xi}(n)$ has the spectral representation
 $$\vec{\xi}(n) = \int_{-\pi}^{\pi}e^{in\theta}d\vec{Z}(\theta), n \in \mathbb{Z},$$

\noindent and the covariation has the representation
$$[\vec{\xi}(n),\vec{\xi}(m)]_\alpha = \int_{-\pi}^{\pi}e^{i(n-m)\theta}d \mu(\theta), m, n \in \mathbb{Z}.$$
\end{definition}

\begin{definition} [periodically harmonizable $S\alpha S$ sequence]
  A symmetric $\alpha$-stable, $S\alpha S$, stochastic sequence $\{\xi(n),n\in\mathbb Z\}$ is said to be periodically harmonizable, $PHS{\alpha}S$, if
  the vector-valued stochastic sequence \linebreak
 $\vec{\xi}(n)=\left \{ \xi(nT+k) \right \}_{k = 0} ^ {T-1},\,n\in \mathbb{Z}$,
 is harmonizable symmetric $\alpha$-stable, $HS{\alpha}S$, stochastic sequence.
\end{definition}

Note that a $HS\alpha S$ stochastic sequence is not necessarily stationary even second order stationary, but for $\alpha = 2$ the $HS\alpha S$ sequences are stationary with Gaussian distribution.

In this article we consider the case where $1<\alpha \leq 2$.

Denote by $H(\xi)$ the time domain of the $HS\alpha S$ sequence $\{\xi(n),n\in\mathbb Z\}$, which is a closed in the norm $\|\cdot\|_{\alpha}$ linear manifold generated by all values of the $HS\alpha S$ sequence $\{\xi(n),n\in\mathbb Z\}$.
It follows from the spectral representation of the $HS\alpha S$ sequence $\{\xi(n),n\in\mathbb Z\}$
that the mapping $\xi(n)\leftrightarrow e^{in\theta},n\in\mathbb Z,$ extents to an isomorphism between the spaces $H(\xi)$ and $L^{\alpha}(\mu)$. Under this isomorphism to each $\eta \in H(\xi)$ corresponds a unique $f\in L^{\alpha}(\mu)$ such that $\eta=\int_{-\pi}^{\pi}f(\theta)dZ(\theta)$.

For a closed linear subspace $M \subseteq L^{\alpha}(\mu)$ and $f \in L^{\alpha}(\mu)$, there exists a unique element from $M$ which minimizes the distance to $f$. This element is called projection of $f$ onto $M$ or the best approximation of $f$ in $M$. This projection is denoted by $P_M f$ and is uniquely determined by the condition (Singer, 1970)
\begin{equation}
 \int_{-\pi}^{\pi} g \left(f - P _M f\right)^{<\alpha - 1>}d \mu = 0,\quad g \in M.
\end{equation}

Similarly, for $HS \alpha S$ stochastic sequence $\{\xi(n),n\in\mathbb Z\}$ and a closed linear subspace $H^-(\xi)$ of the space $H(\xi)$
there is a uniquely determined element $\hat{\xi}(n) \in H^-(\xi)$ which minimizes the distance to $\xi(n)$ and is uniquely determined from the condition
 \begin{equation}\label{eq:ortogonal}
 \left[\eta, \xi(n) - \hat{\xi}(n)\right]_{\alpha} = 0,\quad \eta \in H^-(\xi).
 \end{equation}

From linearity of the covariation with respect to the first argument from this relation we have that
 \begin{multline}\label{eq:linearity}
 ||\xi(n) - \hat{\xi}(n)||_{\alpha}^{\alpha} = \left[\xi(n),\xi(n)-\hat{\xi}(n)\right]_{\alpha}-\left[\hat{\xi}(n), \xi(n)-\hat{\xi}(n)\right]_{\alpha} =
  \\
  =\left[\xi(n), \xi(n)-\hat{\xi}(n)\right]_{\alpha}.
 \end{multline}
This relation plays a fundamental role in the characterization of minimal $HS \alpha S$ stochastic sequences $\{\xi(n),n\in\mathbb Z\}$.

\section{Interpolation problem. Observations with noise. Projection approach}

Consider the problem of the optimal estimation of the linear functional
$$A_N \vec{\xi} =\sum_{j = 0}^{N} (\vec{a}(j))^{\top} \vec{\xi}(j) = \int_{-\pi}^{\pi} (A_N(e^{i\theta}))^{\top}d \vec{Z}^{\xi}(\theta),$$
where
$$A_N(e^{i\theta})=\sum_{j = 0}^{N}\vec{a}(j) e^{ij\theta},
$$
that depends on the
unknown values  $\vec{\xi}(j),j=0,1,\dots,N,$  of a vector-valued harmonizable symmetric $\alpha$-stable random sequence $\vec{\xi}(j)=\left \{ \xi_ {k} (j) \right \}_{k = 1} ^ {T}$,
from observations of the sequence $\vec{\xi}(j)+\vec{\eta}(j)$ at points  $j\in\mathbb Z\setminus\{0,1,\dots,N\}$.

We consider the problem for mutually independent
vector-valued harmonizable symmetric $\alpha$-stable random sequences $\vec{\xi}(j)=\left \{ \xi_ {k} (j) \right \}_{k = 1} ^ {T}$
and
 $\vec{\eta}(j)=\left \{ \xi_ {k} (j) \right \}_{k = 1} ^ {T}$
 which have
absolutely continuous spectral measures and the spectral density $f(\theta)$ and $g(\theta),$
satisfying the minimality condition (Pourahmadi, 1984; Weron, 1985)
   \begin{equation}\label{minimalityF+G}
   \int_{-\pi}^{\pi} \text{Tr}\,\left[(f(\theta)+g(\theta))^{-1/(\alpha-1)}\right]d\theta<\infty.
    \end{equation}

Denote by $H^N(\xi+\eta)$ the closed in the $||\cdot||_\alpha$ norm linear manifold generated by values of the harmonizable symmetric $\alpha$-stable stochastic sequence
$\vec{\xi}(j)+\vec{\eta}(j)$, $j\in\mathbb Z\setminus\{0,1,\dots,N\}$
 in the space $H(\xi+\eta)$ generated by all values of the harmonizable symmetric $\alpha$-stable, $HS \alpha S$, stochastic sequence $\vec{\xi}(j)+\vec{\eta}(j),\,j\in \mathbb{Z}$.

The optimal estimate $\hat{A}_N \vec{\xi}$ of the functional ${A}_N \vec{\xi}$ is a projection of ${A}_N \vec{\xi}$ on the subspace $H^N(\xi+\eta)$
which is determined by the relations
 $$[\eta, A_N \vec{\xi} - \hat{A}_N \vec{\xi}\,]_{\alpha} = 0, \quad \forall \eta \in H^N(\xi+\eta),$$
or, equivalently, by relations
  \begin{equation}\label{ortogon1F+G}
  [\vec{\xi}(j)+\vec{\eta}(j), A_N \vec{\xi} - \hat{A}_N \vec{\xi}]_{\alpha} = 0, \quad j\in\mathbb Z\setminus\{0,1,\dots,N\}.
   \end{equation}
It follows from the isomorphism between the spaces $H(\xi+\eta)$ and $L^{\alpha}(F+G)$
that the optimal estimate $\hat{A}_N \vec{\xi}$ of the functional ${A}_N \vec{\xi}$ is of the form
\begin{equation}\label{estimF+G}
\hat{A}_N \vec{\xi} = \int_{-\pi}^{\pi} ({h}(\theta))^{\top} \left(d \vec{Z}^{\xi}(\theta) + d\vec{Z}^{\eta}(\theta) \right).
 \end{equation}
 It is  determined by the spectral characteristic ${h}(\theta)$ of the estimate which is from the subspace $L_N^{\alpha}(F+G)$ of the $L^{\alpha}(F+G)$ space generated by functions
  $$e^{ij\theta}\delta_{k} , \; \delta_{k} = \left \{ \delta_{kl} \right \}_{l=1}^{T} , \; k= {1,\dots,T},\,\, j\in\mathbb Z\setminus\{0,1,\dots,N\}.$$

 The spectral characteristic $h({\theta})=\left\{h_k({\theta})\right\}_{k=1}^T$ of the optimal estimate satisfies the
 following equations
 \begin{multline}\label{ortogon2F+G}
    \int_{-\pi}^{\pi} {e^{ij\theta } }
        \left[f(\theta)\left( A(e^{i\theta}) - {h}(\theta) \right)^{<\alpha - 1>}
 - g(\theta)\left({h}(\theta) \right)^{<\alpha - 1>}\right]
   d\theta = 0, \\
    j\in\mathbb Z\setminus\{0,1,\dots,N\}.
    \end{multline}

It follows from these equations that the spectral characteristic ${h}(\theta)$ of the estimate is determined by the equation
   \begin{equation}\label{spcheq1F+G}
  \left[f(\theta)\left( A(e^{i\theta}) - {h}(\theta) \right)^{<\alpha - 1>}
 - g(\theta)\left({h}(\theta) \right)^{<\alpha - 1>}\right]  = C_N(e^{i\theta}),
 \end{equation}
 $$ C_N(e^{i\theta})=\sum_{j = 0}^{N} \vec{c}(j) e^{-i j \theta},$$
    \noindent where
 $\vec{c}(j)=\left \{{c}_k(j)\right \}_{k=1}^{T} ,j=0,1,\dots,N$ are unknown coefficients.

The unknown coefficients $\vec{c}(j),j=0,1,\dots,N,$ are determined from the condition
${h}(\theta)\in L_N^{\alpha}(f+g)$,
which gives us the system of equations
\begin{equation}\label{sp_eq1F+G}
\int_{-\pi}^{\pi} e^{-i\theta k}\,\, {h}(\theta) d\theta = 0,\quad  k = 0, 1,\dots, N.
\end{equation}

 The variance of the optimal estimate of the functional is calculated by the formula
\begin{multline}\label{var2F+G}
\left\|{A}_N \xi- \hat{A}_N \xi \right\|_\alpha^\alpha =
\\
=
\int_{-\pi}^{\pi}
\left( A_N(e^{i\theta}) - {h}(\theta) \right)^{\top}
f(\theta)
\left( A_N(e^{i\theta}) - {h}(\theta) \right)^{<{\alpha - 1}>}
d \theta+
\\
 +
 \int_{-\pi}^{\pi}
\left( {h}(\theta) \right)^{\top}
g(\theta)
\left(  {h}(\theta) \right)^{<{\alpha - 1}>}
d \theta
 \end{multline}

In the case where the dimension of the sequence $T=1$, this formula can be written in the form
\begin{equation}\label{T=1var2F+G}
\left\|{A}_N \xi- \hat{A}_N \xi \right\|_\alpha^\alpha =
\int_{-\pi}^{\pi}
\left| A_N(e^{i\theta}) - {h}(\theta) \right|^{\alpha}
f(\theta)
d \theta+
 \int_{-\pi}^{\pi}
\left| {h}(\theta) \right|^{\alpha}
g(\theta)
d \theta.
 \end{equation}

We can conclude that the following theorem holds true.

\begin{thm}\label{thm11}
Let $\vec{\xi}(j)=\left \{ \xi_ {k} (j) \right \}_{k = 1} ^ {T},\,j\in \mathbb{Z}$  and $\vec{\eta}(j)=\left \{ \eta_ {k} (j) \right \}_{k = 1} ^ {T},\,j\in \mathbb{Z}$
be mutually independent harmonizable symmetric $\alpha$-stable $HS{\alpha}S$, stochastic sequences which have
absolutely continuous spectral measures and the spectral densities $f(\theta)$ and $g(\theta)$ satisfying the minimality condition (\ref{minimalityF+G}).
The optimal linear estimate  $\hat{A}_N \vec{\xi}$ of the functional $A_N \vec{\xi} =\sum_{j = 0}^{N} (\vec{a}(j))^{\top} \vec{\xi}(j)$,
that depends on the unknown values $\vec{\xi}(j),j=0,1,\dots,N,$  of the sequence $\vec{\xi}(j)$, from observations of the sequence
 $\vec{\xi}(j)+\vec{\eta}(j)$ at points  $j\in\mathbb Z\setminus\{0,1,\dots,N\}$
 is calculated by formula (\ref{estimF+G}).
 The spectral characteristic ${h}(\theta)$ of the estimate is determined by  equation
(\ref{spcheq1F+G}), where the unknown coefficients $\vec{c}(j),j=0,1,\dots,N,$ are determined from the system of equations
(\ref{sp_eq1F+G}).
The variance of the optimal estimate of the functional is calculated by formula (\ref{var2F+G}) (by formula (\ref{T=1var2F+G}) in the case where dimension $T=1$).
\end{thm}

\begin{thm}\label{thm22}
Let the vector-valued harmonizable symmetric $\alpha$-stable, $HS{\alpha}S$, stochastic sequences
 $\vec{\xi}(n)=\left \{ \xi(nT+k) \right \}_{k = 0} ^ {T-1}$,
  $\vec{\eta}(n)=\left \{ \eta(nT+k) \right \}_{k = 0} ^ {T-1},\,n\in \mathbb{Z}$,
 that
 correspond to
  periodically harmonizable symmetric $\alpha$-stable, $PHS{\alpha}S$, stochastic sequences
  $\{\xi(n),n\in\mathbb Z\}$  and $\{\eta(n),n\in\mathbb Z\}$,
  have
absolutely continuous spectral measures and the spectral densities $f(\theta),$ $g(\theta),$ satisfying the minimality condition (\ref{minimalityF+G}).
The optimal linear estimate  $\hat{A}_N \vec{\xi}$ of the functional
$$A_N \vec{\xi} =\sum_{j = 0}^{N} (\vec{a}(j))^{\top} \vec{\xi}(j)=\sum_{j = 0}^{N} \sum_{k = 0}^{T-1}
{a}_k(j) {\xi}(jT+k),$$
where $\vec{a}(j)= \{{a}_k(j)\}_{k = 0}^{T-1},$ $\vec{\xi}(j)=\{{\xi}(jT+k)\}_{k = 0}^{T-1}$,
that depends on the unknown values ${\xi}(j),j=0,1,\dots,T(N+1)-1,$ of the sequence ${\xi}(j)$ from observations of the sequence ${\xi}(j)+{\eta}(j)$
at points $j\in\mathbb Z\setminus\{0,1,\dots,T(N+1)-1\}$
is calculated by the formula (\ref{estimF+G}).
 The spectral characteristic ${h}(\theta)$ of the estimate is determined by  equation
(\ref{spcheq1F+G}), where the unknown coefficients  $\vec{c}(j)= \{{c}_k(j)\}_{k = 0}^{T-1},j=0,1,\dots,T(N+1)-1$, are determined from the system of equations (\ref{sp_eq1F+G}).
The variance of the optimal estimate of the functional is calculated by the formula (\ref{var2F+G})  (by formula (\ref{T=1var2F+G}) in the case where period $T=1$).
\end{thm}

\subsection{Interpolation problem. Observations without noise. Projection approach}

Consider the problem of the optimal estimation of the linear functional
$$A_N \vec{\xi} =\sum_{j = 0}^{N} (\vec{a}(j))^{\top} \vec{\xi}(j) = \int_{-\pi}^{\pi} (A_N(e^{i\theta}))^{\top}d\vec{Z}^{\xi}(\theta),$$
where
$$A_N(e^{i\theta})=\sum_{j = 0}^{N}\vec{a}(j) e^{ij\theta},
$$
that depends on the
unknown values  $\vec{\xi}(j),j=0,1,\dots,N,$  of a vector-valued harmonizable symmetric $\alpha$-stable random sequence $\vec{\xi}(j)=\left \{ \xi_ {k} (j) \right \}_{k = 1} ^ {T}$,  from observations of the sequence $\vec{\xi}(j)$ at points  $j\in\mathbb Z\setminus\{0,1,\dots,N\}$.

We consider the problem for
vector-valued harmonizable symmetric $\alpha$-stable random sequence $\vec{\xi}(j)=\left \{ \xi_ {k} (j) \right \}_{k = 1} ^ {T}$
 which have
absolutely continuous spectral measure $\mu$ and the spectral density $f(\theta),$
satisfying the minimality condition (Pourahmadi, 1984; Weron, 1985)
   \begin{equation}\label{minimality}
   \int_{-\pi}^{\pi} \text{Tr}\,\left[(f(\theta))^{-1/(\alpha-1)}\right]d\theta<\infty.
    \end{equation}

Denote by
$H^N(\xi)$ the closed in the $||\cdot||_\alpha$ norm linear manifold generated by values of the harmonizable symmetric $\alpha$-stable random sequence
$\vec{\xi}(j)$, $j\in\mathbb Z\setminus\{0,1,\dots,N\}$
 in the space $H(\xi)$ generated by all values of the $HS \alpha S$ sequence $\vec{ \xi}(j)=\left \{ \xi_ {k} (j) \right \}_{k = 1} ^ {T},\,j\in \mathbb{Z}$.

The optimal estimate $\hat{A}_N \vec{\xi}$ of the functional ${A}_N \vec{\xi}$ is a projection of ${A}_N \vec{\xi}$ on the subspace $H^N(\xi)$
which is determined by the relations
 $$[\eta, A_N \vec{\xi} - \hat{A}_N \vec{\xi}\,]_{\alpha} = 0, \quad \forall \eta \in H^N(\xi),$$
or, equivalently, by relations
  \begin{equation}\label{ortogon1}
  [\vec{\xi}(j), A_N \vec{\xi} - \hat{A}_N \vec{\xi}]_{\alpha} = 0, \quad j\in\mathbb Z\setminus\{0,1,\dots,N\}.
   \end{equation}
It follows from the isomorphism between the spaces $H(\xi)$ and $L^{\alpha}(\mu)$
that the optimal estimate $\hat{A}_N \vec{\xi}$ of the functional ${A}_N \vec{\xi}$ is of the form
\begin{equation}\label{estim}
\hat{A}_N \vec{\xi} = \int_{-\pi}^{\pi} ({h}(\theta))^{\top} d\vec{Z}^{\xi}(\theta).
 \end{equation}
 It is  determined by the spectral characteristic ${h}(\theta)$ of the estimate which is from the subspace $L_N^{\alpha}(\mu)$ of the $L^{\alpha}(\mu)$ space generated by functions
  $$e^{ij\theta}\delta_{k} , \; \delta_{k} = \left \{ \delta_{kl} \right \}_{l=1}^{T} , \; k= {1,\dots,T}, j\in\mathbb Z\setminus\{0,1,\dots,N\}.$$

 The spectral characteristic $h({\theta})=\left\{h_k({\theta})\right\}_{k=1}^T$ of the optimal estimate satisfies the
 following equations
 \begin{equation}\label{ortogon2}
    \int_{-\pi}^{\pi} {e^{ij\theta } f(\theta)\left( A_N(e^{i\theta}) - {h}(\theta) \right)^{<\alpha - 1>}}d\theta = 0, \quad
    j\in\mathbb Z\setminus\{0,1,\dots,N\}.
    \end{equation}

It follows from these equations that the spectral characteristic ${h}(\theta)$ of the estimate is determined by the equation
    $$f(\theta)\left( A_N(e^{i\theta}) - {h}(\theta) \right)^{<\alpha - 1>}  = C_N(e^{i\theta}),\quad C_N(e^{i\theta})=\sum_{j = 0}^{N} \vec{c}(j) e^{-i j \theta},$$
    \noindent where
 $\vec{c}(j)=\left \{{c}_k(j)\right \}_{k=1}^{T} ,j=0,1,\dots,N$ are unknown coefficients.

 It follows from the last relation that the spectral characteristic ${h}(\theta)$ of the estimate is of the form
    \begin{equation}\label{spectralcharacteristic}
    {h}(\theta) = A_N(e^{i\theta}) - \left(f^{-1}(\theta)C_N(e^{i\theta}) \right)^{<\frac{1}{\alpha - 1}>}.
    \end{equation}
The unknown coefficients $\vec{c}(j),j=0,1,\dots,N,$ are determined from the condition
${h}(\theta)\in L_N^{\alpha}(\mu)$,
which gives us the system of equations
\begin{equation}\label{sp_eq1}
\int_{-\pi}^{\pi} e^{-i\theta k}\,\, {h}(\theta) d\theta = 0,\quad  k = 0, 1,\dots, N.
\end{equation}
These equations are of the form
\begin{multline}\label{sp_eq2}
  \int_{-\pi}^{\pi} e^{-i\theta k} \left[ \left(\sum_{j = 0}^{N}\vec{a}(j) e^{i j \theta}\right) - \left(\left( f(\theta) \right)^{-1}\left(\sum_{j = 0}^{N} \vec{c}(j) e^{-i j \theta}\right) \right)^{<\frac{1}{\alpha - 1}>}\right]   d\theta = 0,
  \\
   k = 0, 1,\dots, N.
  \end{multline}

 The variance of the optimal estimate of the functional is calculated by the formula

    \begin{multline}\label{variance2}
\left\| {A}_N \vec{\xi}- \hat{A}_N \vec{\xi} \right\|_\alpha^\alpha =
\\
=
\int_{-\pi}^{\pi}
\left[\left(f^{-1}(\theta)C_N(e^{i\theta}) \right)^{<\frac{1}{\alpha - 1}>}\right]^{\top}
f(\theta)
\left(f^{-1}(\theta)C_N(e^{i\theta}) \right)^{<\frac{\alpha-1}{\alpha - 1}>}
d \theta.
\end{multline}

In the case where the dimension of the sequence $T=1$, this formula can be written in the form
\begin{equation}\label{T=1variance2}
\left\|{A}_N \xi- \hat{A}_N \xi \right\|_\alpha^\alpha =
\int_{-\pi}^{\pi}
\left|
\left(f^{-1}(\theta)C_N(e^{i\theta}) \right)^{<\frac{1}{\alpha - 1}>}
\right|^{\alpha}
f(\theta)
d \theta.
 \end{equation}

 We can conclude that the following theorem holds true.

\begin{thm}\label{thm1}
Let $\vec{\xi}(j)=\left \{ \xi_ {k} (j) \right \}_{k = 1} ^ {T},\,j\in \mathbb{Z}$  be a vector-valued harmonizable symmetric $\alpha$-stable, $HS{\alpha}S$,  stochastic sequence
which has
absolutely continuous spectral measure $\mu(\theta)$ and the spectral density $f(\theta),$ satisfying the minimality condition (\ref{minimality}).
The optimal linear estimate  $\hat{A}_N \vec{\xi}$ of the functional $A_N \vec{\xi} =\sum_{j = 0}^{N} (\vec{a}(j))^{\top} \vec{\xi}(j)$,
that depends on the unknown values $\vec{\xi}(j),j=0,1,\dots,N,$  of the sequence $\vec{\xi}(j)$, from observations of the sequence $\vec{\xi}(j)$
at points $j\in\mathbb Z\setminus\{0,1,\dots,N\}$ is calculated by the formula
(\ref{estim}). The spectral characteristic ${h}(\theta)$ of the optimal estimate is determined by the equation
(\ref{spectralcharacteristic}), where the unknown coefficients $\vec{c}(j),j=0,\dots, N,$ are determined from the system of equations
(\ref{sp_eq2}).
The variance of the optimal estimate of the functional is calculated by the formula (\ref{variance2}) (by formula (\ref{T=1variance2}) in the case where dimension of the sequence  $T=1$).
\end{thm}

\begin{thm}\label{thm2 fs}
Let the vector-valued harmonizable symmetric $\alpha$-stable, $HS{\alpha}S$, stochastic sequence
 $\vec{\xi}(n)=\left \{ \xi(nT+k) \right \}_{k = 0} ^ {T-1},\,n\in \mathbb{Z}$, that corresponds to
  periodically harmonizable symmetric $\alpha$-stable, $PHS{\alpha}S$, stochastic sequence
  $\{\xi(n),n\in\mathbb Z\}$,
 has
absolutely continuous spectral measure $\mu(\theta)$ and the spectral density $f(\theta),$ satisfying the minimality condition (\ref{minimality}).
The optimal linear estimate  $\hat{A}_N \vec{\xi}$ of the functional
$$A_N \vec{\xi} =\sum_{j = 0}^{N} (\vec{a}(j))^{\top} \vec{\xi}(j)=\sum_{j = 0}^{N} \sum_{k = 0}^{T-1}
{a}_k(j) {\xi}(jT+k),$$
 $\vec{a}(j)= \{{a}_k(j)\}_{k = 0}^{T-1},$ $\vec{\xi}(j)=\{{\xi}(jT+k)\}_{k = 0}^{T-1}$,
that depends on the unknown values ${\xi}(j),j=0,1,\dots,T(N+1)-1,$ of the sequence ${\xi}(j)$ from observations of the sequence ${\xi}(j)$
at points $j\in\mathbb Z\setminus\{0,1,\dots,T(N+1)-1\}$
is calculated by the formula (\ref{estim}).
The spectral characteristic
${h}(\theta)$ of the optimal estimate is of the form
(\ref{spectralcharacteristic}), where the unknown coefficients $\vec{c}(j)= \{{c}_k(j)\}_{k = 0}^{T-1},j=0,1,\dots,T(N+1)-1$, are determined from the system of equations (\ref{sp_eq2}).
The variance of the optimal estimate of the functional is calculated by the formula (\ref{variance2}) (by formula (\ref{T=1variance2}) in the case where period $T=1$).
\end{thm}

\begin{exm}
Consider the problem of the optimal linear estimation of the functional
$A_0 \xi = a \xi(0)$,
that depends on the unknown value $\xi(0)$ of a harmonizable symmetric $\alpha$-stable, $HS{\alpha}S$, stochastic sequence $\xi(n)$
that has
absolutely continuous spectral measure $\mu(\theta)$ and the spectral density $f(\theta),$ satisfying the minimality condition (\ref{minimality}),
from observations of the sequence $\xi(n)$  at points $n\in\mathbb Z\setminus\{0\}$ (Moklyachuk, Ostapenko, 2015, 2016).

In this case the spectral characteristic
${h}(\theta)$ of the optimal estimate of the functional is of the form
 $${h}(\theta) = a - c^{<\frac{1}{\alpha - 1}>}\left(f(\theta)\right)^\frac{-1}{\alpha - 1}.$$
The variance of the optimal estimate of the functional is calculated by the formula
$$\left\| \hat{A}_0 \xi - A_0 \xi \right\|_\alpha^\alpha = \int_{-\pi}^{\pi} \left|c^{<\frac{1}{\alpha - 1}>} \left(f(\theta)\right)^{-\frac{1}{\alpha - 1}}\right|^{\alpha} f(\theta) d \theta,$$
where the constant $c$ is  a solution of the equation
$$\int_{-\pi}^{\pi} \left(a - c^{<\frac{1}{\alpha - 1}>} \left(f(\theta)\right)^{-\frac{1}{\alpha - 1}}\right) d \theta = 0,$$
$$c = \frac{\left(2 \pi a \right)^{<\alpha - 1>}}{\left(\int_{-\pi}^{\pi}{(f(\theta))^{-\frac{1}{\alpha - 1}}} d \theta \right)^{<\alpha - 1>} } .$$

In the case of $\alpha=2$ we have the Kolmogorov (see selected works by Kolmogorov (1992)) results
$${h}(\theta) = a - c\left(f(\theta)\right)^{-1},$$
$$c = a \left(\frac{1}{2\pi}\int_{-\pi}^{\pi}{\left((f(\theta))^{-1}\right)} d \theta \right)^{-1},$$
$$\left\| \hat{A} \xi - A \xi \right\|_2^2 =
2\pi |a|^{2} \left(\frac{1}{2\pi}\int_{-\pi}^{\pi}{\left((f(\theta))^{-1}\right)} d \theta \right)^{-1}.$$
\end{exm}

\begin{exm}
Consider the problem of the optimal  linear estimation of the functional
$A_1 \xi = a(0)\xi(0) + a(1)\xi(1)$, that depends on the unknown values $\xi(0)$ and  $\xi(1)$  of a harmonizable symmetric $\alpha$-stable, $HS{\alpha}S$, stochastic sequence
$\xi(n), n \in \mathbb{Z}$ with $\alpha = \frac{4}{3}$  and the spectral density $f(\theta) = |e^{i\theta} + d|^{-\frac{4}{3}}, -1 < d < 1,$
from observations of the sequence $\xi(n)$  at points $n\in\mathbb Z\setminus\{0,1\}$ (Moklyachuk, Ostapenko, 2015, 2016).

In this case the spectral characteristic
${h}(\theta)$ of the optimal estimate of the functional is of the form
\begin{equation}\label{eq:spectral_chaacteristique_example}
{h}(\theta) = a_0 + a_1e^{i \theta} -  \left(c_0 + c_1 e^{-i \theta}\right)^{<3>} |e^{i\theta} + d|^{4}.
\end{equation}
Here
\begin{multline}
\label{eq:coeff_c1}\left(c_0 + c_1 e^{-i \theta}\right)^{<3>} =
 \\
 =\left(c_0 + c_1 e^{-i \theta}\right) \left(\overline{c_0} + \overline{c_1} e^{i \theta}\right)^2 = b_{-1} e^{-i \theta} + b_0 + b_{1} e^{i \theta} + b_{2} e^{2 i \theta},
\end{multline}
where
\begin{equation}\label{eq:coeff_c2}
b_{-1}= \overline{c_0}^2 c_1,\;b_{0} = c_0 \overline{c_0}^2 + 2\overline{c_0}c_1\overline{c_1},\;b_1 = 2 c_0 \overline{c_0}\overline{c_1} + c_1 \overline{c_1}^2,\;b_2 = c_0 \overline{c_1}^2,
\end{equation}
and
\begin{equation}\label{eq:spectral_density_example11}
|e^{i \theta} + d|^{4} = r_{-2} e^{- 2 i\theta} + r_{-1} e^{- i\theta} + r_{0} + r_{1} e^{i\theta} + r_{2} e^{2 i\theta},
\end{equation}
where
\begin{equation}\label{eq:spectral_density_example22}
r_{-2} = d^2,\;r_{-1} =2 d + 2 d^3,\;r_0 =1 + 4 d^2 + d^4,\;r_1 = 2 d + 2 d^3,\;r_2 = d^2.
\end{equation}
It follows from (\ref{eq:coeff_c1}) -- (\ref{eq:spectral_density_example22}) that the spectral characteristic
${h}(\theta)$ of the optimal estimate of the functional is of the form
$${h}(\theta) = h_{-3} e^{-3 i \theta} + h_{-2} e^{- 2 i \theta} + h_{-1} e^{-i \theta} + h_{0} + h_{1} e^{i \theta} + h_{2} e^{2 i \theta} + h_3 e^{3 i \theta} + h_{4} e^{4 i \theta},$$
where
$$h_{-3} = -b_{-1} r_{-2},\;h_{-2} = -b_{-1} r_{-1} -b_{0} r_{-2},\;h_{-1} = -b_{-1} r_{0} - b_{0} r_{-1} - b_{1} r_{-2},$$
 $$h_0 = a_0 - b_{-1} r_{1} - b_{0} r_{0} - b_{1} r_{-1} - b_{2} r_{-2},\;h_1 = a_1 - b_{-1} r_{2} - b_{0} r_{1} - b_{1} r_{0}- b_{2} r_{-1},$$
 $$h_{2} =  -b_{0} r_{2} - b_{1} r_{1} - b_{1} r_{0},\;h_3 = -b_{1}r_{2} - b_{2} r_{1},\;h_4 = -b_{2}r_{2}.$$
 Condition (\ref{sp_eq1}) is satisfied if
\begin{multline}\label{eq:coeff_zero}
h_0 = a_0 - b_{-1} r_{1} - b_{0} r_{0} - b_{1} r_{-1} - b_{2} r_{-2}= 0,\\
h_1= a_1 - b_{-1} r_{2} - b_{0} r_{1} - b_{1} r_{0}- b_{2} r_{-1} = 0.
\end{multline}
These equations and equations (\ref{eq:coeff_c2}) and (\ref{eq:spectral_density_example22}) give us a system of equations that determine the unknown coefficients $c_0, c_1.$

Consider the problem for  $a(0) = 1, a(1) = 1, d = 0.5$.
In this case the unknown coefficients are calculated from the indicated equations. They are as follows: $ c_0 \approx 0.44,$ $ c_1 \approx 0.44.$
The spectral characteristic of the optimal estimate of the functional is of the form
$$
{h}(\theta) = h_{-3} e^{-3 i \theta} + h_{-2} e^{- 2 i \theta} + h_{-1} e^{-i \theta} + h_{2} e^{2 i \theta} + h_3 e^{3 i \theta} + h_{4} e^{4 i \theta},
$$ with
\begin{multline*}
h_{-3} \approx -0.02,~h_{-2} \approx -0.17,~h_{-1} \approx -0.57,
\\
~h_{2} \approx -0.57,~h_{3} \approx -0.17,~h_{4} \approx -0.02.
\end{multline*}
The variation  of the optimal estimate of the functional is calculated by the formula
\begin{multline*}
\left\| \hat{A}_1 \xi - A_1 \xi \right\|_{\frac{4}{3}}^{\frac{4}{3}} \approx
\\
\int_{-\pi}^{\pi} \left|\left(0.44 +0.44e^{i \theta}\right)^{<3>} \left| e^{i\theta} + 0.5 \right|^4\right|^{\frac{4}{3}} \left(|e^{i\theta} + 0.5|^{-\frac{4}{3}}\right) d \theta \approx 5.57.
\end{multline*}

The corresponding results for $\alpha = 2$ and $f(\theta) = |e^{i\theta} + 0.5|^{-2}$ are the following: $ c_0 = \frac{4}{7},$ $ c_1 = \frac{4}{7},$
\begin{multline*}
\left\| \hat{A}_1 \xi - A_1 \xi \right\|_{2}^{2} =
 \\
 \int_{-\pi}^{\pi} \left|\left(\frac{4}{7} + \frac{4}{7}e^{-i \theta}\right)^{<1>} \left| e^{i\theta} + 0.5 \right|^{2}\right|^{2} \left(|e^{i\theta} + 0.5|^{-2}\right) d \theta  = \frac{16\pi}{7}.
\end{multline*}
\end{exm}

\begin{exm}
Consider the problem of the optimal  linear estimation of the functional
$A_N \xi~=~\sum_{j = 0}^{N} a(j) \xi(j)$, that depends on the unknown values $\xi(j),j=0,1,\dots,N,$
of a harmonizable symmetric $\alpha$-stable, $HS{\alpha}S$, stochastic sequence
$\xi(n), n \in \mathbb{Z}$ with $\alpha = \frac{4}{3}$  and the spectral density $f(\theta)$
from observations of the sequence $\xi(n)$  at points $n\in\mathbb Z\setminus\{0,1,\dots,N\}$ (Moklyachuk, Ostapenko, 2015, 2016).

In this case the spectral characteristic
${h}(\theta)$ of the optimal estimate of the functional is of the form
\begin{equation}\label{eq:spectral_chaacteristique_example}
{h}(\theta) = \sum_{k = 0}^N a_k e^{i k \theta} -  \left(\sum_{k = 0}^N c_k e^{-i k \theta}\right)^{<3>}  (f(\theta))^{-3}.
\end{equation}
Making use of the Fourier representation
$$\left(\sum_{k = 0}^N c_k e^{-i k \theta}\right)^{<3>} = \left(\sum_{k = 0}^N c_k e^{-i k \theta}\right) \left(\sum_{k = 0}^N \overline{c}_k e^{ i k \theta}\right)^2 = \sum_{k=-N}^{2N}b_k e^{i k \theta}$$
and
$$(f(\theta))^{-3}=\sum_{k=-\infty}^{\infty}r_k e^{i k \theta}, \quad r_k = \frac{1}{2 \pi} \int_{-\pi}^{\pi}(f(\theta))^{-3}e^{-ik\theta} d\theta,$$
we have
\begin{equation}\label{eq:spectr_example}
{h}(\theta) = \sum_{k = 0}^N a_k e^{i k \theta} -  \left(\sum_{k=-N}^{2N}b_k e^{i k \theta}\right)\left(\sum_{k=-\infty}^{\infty}r_k e^{i k \theta}\right).
\end{equation}

To determine a system of equations for calculation the unknown coefficients $c_0, c_1, \ldots, c_{N}$
we define a $(N +1)\times (2N+1)$ matrix $\bold C_N$ with elements
$C_{k,j}=c_{N+k-j}$ for $k\leq j\leq N+k$ and $C_{k,j}=0$ for $k>j$ and $j>N+k$, where $k=0,1,\dots,N,$ $j=0,1,\dots, 2N$.

As a result of multiplication $\vec{\alpha}=\vec{c}\,\bold C_N$ of the vector
 $\vec{c}=(\overline{c}_0, \overline{c}_1, \ldots, \overline{c}_{N})$ and the matrix $\bold C_N$ we get a $2N+1$ vector
 $\vec{\alpha}=({\alpha}_{-N},\dots,\alpha_{N})$.
 Similarly to how the $\bold C_N$ matrix was defined, we define
  $(N +1)\times (3N+1)$ matrix $\Lambda_N$ with elements $\Lambda_{k,j}=\alpha_{-N-k+j}$
for $k\leq j\leq 2N+k$ and $\Lambda_{k,j}=0$ for $k>j$ and $j>2N+k$, where $k=0,1,\dots,N,$ $j=0,1,\dots, 3N$.
As a result of multiplication  $\vec{b}=\vec{c}\,\Lambda_N$ of the vector $\vec{c}$ and the matrix $\Lambda_N$ we get $3N+1$ vector $\vec{b}=({b}_{-N},\dots,b_{2N})$.

Another equation we will get  (\ref{sp_eq1}) by using the form of the spectral characteristic (\ref{eq:spectr_example}). We will have
$\vec{a} = \vec{b}\, \bold R,$
where
 $\vec{a} = (a_{0}, a_{1}, \ldots, a_{N})$, $\bold R$ is the $(3N+1)\times (N+1)$ matrix with elements
 $\{R_{k,j}\}_{k = 0, j = 0}^{3N, N}$, $R_{k,j} = r_{N-k+j}$, constructed from the Fourier coefficients
$$r_k = \frac{1}{2 \pi} \int_{-\pi}^{\pi}(f(\theta))^{-3}e^{-ik\theta} d\theta$$ of the function $(f(\theta))^{-3}$.

The introduced equations
$$
\vec{\alpha}=\vec{c}\,\bold C_N,\quad
\vec{b}=\vec{c}\,\Lambda_N,\quad
\vec{a} = \vec{b} \,\bold R
$$
determine the unknown coefficients $c_0, c_1, \ldots, c_{N}$.

In the particular case $N = 1$ and $f(\theta) = |e^{i\theta} + d|^{-\frac{4}{3}}$ we will have
$$
{h}(\theta) = \sum_{k = 0}^1 a_k e^{i k \theta} -  \left(\sum_{k=-1}^{2}b_k e^{i k \theta}\right)\left(\sum_{k=-2}^{2}r_k e^{i k \theta}\right)
$$
and equations
$$
b_{-1}= \overline{c_0}^2 c_1,\;b_{0} = c_0 \overline{c_0}^2 + 2\overline{c_0}c_1\overline{c_1},\;b_1 = 2 c_0 \overline{c_0}\overline{c_1} + c_1 \overline{c_1}^2,\;b_2 = c_0 \overline{c_1}^2
$$
$$
 a_0 - b_{-1} r_{1} - b_{0} r_{0} - b_{1} r_{-1} - b_{2} r_{-2}= 0,\; a_1 - b_{-1} r_{2} - b_{0} r_{1} - b_{1} r_{0}- b_{2} r_{-1} = 0,
$$
which coincide with the equations  (\ref{eq:coeff_c2}), (\ref{eq:coeff_zero}).
\end{exm}

\subsection{Interpolation problem. Stationary sequences}

Consider the problem of the optimal estimation of the linear functional
$$A_N \vec{\xi} =\sum_{j = 0}^{N} (\vec{a}(j))^{\top} \vec{\xi}(j) = \int_{-\pi}^{\pi} (A_N(e^{i\theta}))^{\top}d\vec{Z}^{\xi}(\theta),$$
where
$$A_N(e^{i\theta})=\sum_{j = 0}^{N}\vec{a}(j) e^{ij\theta},
$$
that depends on the
unknown values  $\vec{\xi}(j),j=0,1,\dots,N,$  of a vector-valued harmonizable symmetric $\alpha$-stable random sequence $\vec{\xi}(j)=\left \{ \xi_ {k} (j) \right \}_{k = 1} ^ {T}$,  from observations of the sequence $\vec{\xi}(j)$ at points  $j\in\mathbb Z\setminus\{0,1,\dots,N\}$  in the particular case where $\alpha=2$.

In this case the harmonizable symmetric $\alpha$-stable stochastic sequences $\vec{\xi}(j)=\left \{ \xi_ {k} (j) \right \}_{k = 1} ^ {T},j\in\mathbb Z\}$ and $\vec{\eta}(j)=\left \{ \eta_ {k} (j) \right \}_{k = 1} ^ {T},j\in\mathbb Z\}$
are stationary sequences and we have the problem of the optimal estimation of the linear functional
$A_N \vec{\xi} =\sum_{j = 0}^{N} (\vec{a}(j))^{\top} \vec{\xi}(j)$
that depends on the unknown values of a stationary random sequence
from observations of the stationary sequence $\vec{\xi}(j)+\vec{\eta}(j)$ at points $j\in\mathbb Z\setminus\{0,1,\dots,N\}$.

We will suppose that  $ \vec{ \xi}(j)= \left \{ \xi_{k} (j) \right
\}_{k=1}^{T} $ and $ \vec{ \eta}(j)= \left \{ \eta_{k} (j) \right
\}_{k=1}^{T} $
are uncorrelated vector-valued stationary stochastic
sequences with spectral density matrices $f( \lambda )= \left
\{f_{ij} ( \lambda ) \right \}_{i,j=1}^{T} $ and $g( \lambda )=
\left \{g_{ij} ( \lambda ) \right \}_{i,j=1}^{T} $
 satisfying the minimality condition

\begin{equation}\label{int-minimality-stat}
 \int_{- \pi}^{ \pi} \, {\mathrm{Tr}}\,  \left[ \left(f( \lambda )+g( \lambda ) \right)^{-1} \right] \, d \lambda \, < \, \infty .
 \end{equation}

\noindent Denote by $L_{2}^{N-}(f+g)$ the subspace of the space $L_{2}(f+g)$ generated by functions of the form
 \[e^{in \lambda} \delta_{k} , \; \delta_{k} = \left \{ \delta_{kl} \right \}_{l=1}^{T} , \; k= {1,\dots,T}, \; n \in Z \backslash \left \{0, \; \ldots , \; N \right \}. \]
\noindent Every linear estimate $ \hat{A}_N \vec{ \xi}$ \label{linest} of the functional $A_N \vec{ \xi}$ based on observations of the sequence $ \vec { \xi} (j) + \vec { \eta} (j) $
for $j \in Z \backslash \left \{0, \; \ldots , \; N \right \}$
is of the form
 \begin{multline}\label{est-stat}
 \hat{A}_{N} \vec{ \xi}= \int_{- \pi}^{ \pi} h(e^{i \lambda} )^{ \top} (Z^{ \xi} (d \lambda )+Z^{ \eta} (d \lambda ))=
\\
 =\int_{- \pi}^{ \pi} \sum_{k=1}^{T}h_{k} (e^{i \lambda} )(Z_{k}^{ \xi} (d \lambda )+Z_{k}^{ \eta} (d \lambda )),
  \end{multline}

\noindent where $Z^{ \xi} ( \Delta )= \left \{Z_{k}^{ \xi} ( \Delta
) \right \}_{k=1}^{T}$ and $Z^{ \eta} ( \Delta )= \left \{Z_{k}^{
\eta} ( \Delta ) \right \}_{k=1}^{T} $ are orthogonal random
measures \label{vypmira} of sequences $ \vec{ \xi}(j) \, $ and $
\vec{ \eta}(j) \,$ respectively, $h(e^{i \lambda} )= \left \{h_{k}
(e^{i \lambda} ) \right \}_{k=1}^{T}$ is the spectral characteristic
of the estimate $ \hat{A}_{N} \vec{ \xi}$. The function $h(e^{i
\lambda} )$ belongs to $L_{2}^{N-} (f+g)$.

\noindent The mean-square error \label{error} $ \Delta (h;f,g)$ of the linear estimate $ \hat{A}_{N} \vec{ \xi}$ is calculated by the formula
 \begin{multline*}
 \Delta (h;f,g)=E \left|A_{N} \vec{ \xi}- \hat{A}_{N} \vec{ \xi} \right|^{2}= \\
 \frac{1}{2 \pi} \int_{- \pi}^{ \pi}(A_N (e^{i \lambda} )-h(e^{i \lambda}))^{ \top} f( \lambda ) \overline{(A_N (e^{i \lambda} )-h(e^{i \lambda} ))}d \lambda
 \\
+ \frac{1}{2 \pi} \int_{- \pi}^{ \pi} \, h(e^{i \lambda} )^{ \top} g( \lambda ) \overline{h(e^{i \lambda} )}d \lambda,
 \end{multline*}
\noindent where
 \[A_{N} (e^{i \lambda} )= \sum_{j=0}^{N} \vec{a}(j)e^{ij \lambda}  . \]

\noindent The spectral characteristic $h(f,g)$ of the optimal linear
estimate $A_{N} \vec{ \xi}$ minimizes the value of the mean-square
error
 \begin{equation} \label{int-Delta-extrF+G}
 \Delta (f,g)= \Delta (h(f,g);F,G)= \mathop{ \min} \limits_{h \in L_{2}^{N-} (f+g)} \Delta (h;f,g)= \mathop{ \min} \limits_{ \hat{A}_{N} \vec{ \xi}} E \left|A_{N} \vec{ \xi}- \hat{A}_{N} \vec{ \xi} \right|^{2} .
 \end{equation}

\noindent The optimal linear estimate $ \hat{A}_N \vec{ \xi}$ is a
solution to the extremum problem (\ref{int-Delta-extrF+G}). It is
determined by two conditions (Kolmogorov, 1992):

1)$\; \hat{A}_N \vec{ \xi} \in H \left[ \xi_{k} (n)+ \eta_{k} (n), k= \overline{1,T}, n \in Z \backslash \left \{0, \ldots, N \right \}\right];$

2)$\; A_N \vec{ \xi}- \hat{A}_N \vec{ \xi} \perp H \left[ \xi_{k} (n)+ \eta_{k} (n), k= \overline{1,T}, n \in Z \backslash \left \{0, \ldots, N \right \}\right].$

 \noindent Here $H \left[ \xi_{k} (n)+ \eta_{k} (n), k=
\overline{1,T}, n \in Z \backslash \left \{0, \ldots, N \right \}
\right]$ denotes the subspace generated by random variables $[
\xi_{k} (n)+ \eta_{k} (n)$, $k= {1,\dots,T}$, $n \in Z \backslash \left \{0, \ldots, N \right \}]$ in the Hilbert space $L_{2}$ of
random variables with finite second moment and zero first moment.
These conditions help us to find the spectral characteristic
$h(f,g)$ and the mean-square error $ \Delta(f,g)$ of the optimal
linear estimate $ \hat{A}_N \vec{ \xi}$ of the functional $A_N\vec{ \xi}$  in the case where the spectral density matrices
$f(\lambda )$ and $g( \lambda )$ are given and condition
(\ref{int-minimality-stat}) is satisfied.

 The second condition is fulfilled if
 \[
 E \left[ ({A}_N \vec{ \xi}- \hat{A}_N \vec{ \xi}) \overline{ \xi_{k} (j)+\eta_k(j)}\right]=0,\quad   j \in Z \backslash \left \{0, \ldots ,N \right \}, \,\,k= 1,2,\dots,T, \]
\noindent that is
\begin{multline*}
E \left[\int _{-\pi }^{\pi }
\left( (A_N (e^{i\lambda } ))^{\top} Z^{\xi } (d\lambda )-(h(e^{i\lambda } ))^{\top} (Z_{k}^{ \xi} (d \lambda )+Z_{k}^{ \eta} (d \lambda))\right)\right.
\times
\\
\left.\times \int _{-\pi }^{\pi }e^{-ij\lambda } \overline{(Z_{k}^{ \xi} (d \lambda )+Z_{k}^{ \eta} (d \lambda))}  \right]=0,
\end{multline*}
\[j\in Z\backslash \left\{0,\ldots ,N\right\}.\]

  \noindent These equations imply the following equations
 \[ \int_{- \pi}^{ \pi}\left((A_N (e^{i \lambda} )^{ \top} f( \lambda )-h(e^{i \lambda} )^{ \top} (f( \lambda )+g( \lambda ))\right)e^{-ij \lambda}d \lambda=0,
 \]
\[j\in Z\backslash \left\{0,\ldots ,N\right\}.\]
From these equations we have the equation
 \[A_N (e^{i \lambda} )^{ \top} f( \lambda )-h(e^{i \lambda} )^{ \top} (f( \lambda )+g( \lambda ))=C_N (e^{i \lambda} )^{ \top} , \]
\noindent where $C_N (e^{i \lambda} )= \sum_{j=0}^{N} \vec{c}(j)e^{ij \lambda}$, $ \vec{c} = \left(\vec{c}(j),j=0,1,2\dots,N \right) $ are unknown coefficients.

\noindent The last equation gives us a possibility to find a form of the spectral characteristic
of the estimate
\begin{multline} \label{int-sp-char-stat}
h(e^{i \lambda} )^{ \top} =(A_N (e^{i \lambda} )^{ \top} f( \lambda )-C_N (e^{i \lambda} )^{ \top} )(f( \lambda )+g( \lambda ))^{-1} =
\\
=A_N (e^{i \lambda} )^{ \top} -(A_N (e^{i \lambda} )^{ \top} g( \lambda )+C_N (e^{i \lambda} )^{ \top} )(f( \lambda )+g( \lambda ))^{-1} .
  \end{multline}
The first condition is equivalent to the system of equations
 \[ \int_{- \pi}^{ \pi}h(e^{i \lambda} )e^{-ij \lambda} d \lambda =0, \; j=0,1,2, \dots,N \]
\noindent which can be written in the form
 \begin{multline} \label{int-EquationsB}
  \int_{- \pi}^{ \pi}A_N (e^{i \lambda} )^{ \top} f( \lambda )(f( \lambda )+g( \lambda ))^{-1} e^{-ij \lambda} d \lambda=
  \\
 = \int_{- \pi}^{ \pi}C_N (e^{i \lambda} )^{ \top} (f( \lambda )+g( \lambda ))^{-1} e^{-ij \lambda} d \lambda, \\ j=0,1,2, \dots,N.
  \end{multline}

\noindent Let us determine the Fourier coefficients of the matrix-valued functions
\[[(f( \lambda )+g( \lambda ))^{-1} ]^{ \top},\, [f( \lambda )(f( \lambda )+g( \lambda ))^{-1} ]^{
\top},\, [f( \lambda )(f( \lambda )+g( \lambda ))^{-1} g( \lambda )]^{\top}\]

\noindent as follows $ (k,j=0,1,2,\dots,N)$
 \[B(k,j)= B(k-j)=\frac{1}{2 \pi} \int_{- \pi}^{ \pi} \left[(f( \lambda )+g( \lambda ))^{-1} \right]^{ \top} e^{-i(j-k) \lambda} d \lambda, \]
 \[D(k,j)= D(k-j)= \frac{1}{2 \pi} \int_{- \pi}^{ \pi} \left[f( \lambda )(f( \lambda )+g( \lambda ))^{-1} \right]^{ \top} e^{-i(j-k) \lambda} d \lambda, \]
\[R(k,j)= R(k-j)= \frac{1}{2 \pi} \int_{- \pi}^{ \pi} \left[f( \lambda )(f( \lambda )+g( \lambda ))^{-1} g( \lambda ) \right]^{ \top} e^{-i(j-k) \lambda} d \lambda. \]
\noindent Making use of these Fourier coefficients we find the following system of equations to determine the unknown coefficients
$ \mathbf{c}_N=\vec{c}_N = \left(\vec{c}(j),j=0,1,2\dots,N \right) $
\begin{equation}\label{mult-seq-extr-equation-noise}
\mathbf{D}_N\mathbf{a}_N=\mathbf{B}_N\mathbf{c}_N,
\end{equation}
Solution to this system is of the form
\begin{equation}
\mathbf{c}_N=\mathbf{B}_N^{-1}\mathbf{D}_N\mathbf{a}_N,
\end{equation}
that is the components of the vector $ \mathbf{c}_N=\vec{c}_N $ are calculated by the formula
\begin{equation}
\vec{c}(j)=(\mathbf{B}_N^{-1}\mathbf{D}_N\vec{a}_N)(j),\,j=0,1,2\dots,N.
\end{equation}
Here vectors
$$ \mathbf{a}_N=\vec{a}_N = \left(\vec{a}(j),j=0,1,2\dots,N \right) ;$$
$$ \mathbf{c}_N=\vec{c}_N = \left(\vec{c}(j),j=0,1,2\dots,N \right) ;$$
matrices
$$\mathbf{B}_N=\{B(k,j)\}_{k,j=0}^{N},\quad  \mathbf{D}_N=\{D(k,j)\}_{k,j=0}^{N}, \quad \mathbf{R}_N =\{R(k,j)\}_{k,j=0}^{N} $$
are  constructed from the corresponding block-matrices $B(k,j)$, $D(k,j)$, $R(k,j)$ of dimension $T \times T$.

 The spectral characteristic $h(e^{i\lambda})$ and the mean-square error $\Delta(f,g)$ of the optimal linear estimate of the functional $A_N\vec{\xi}$
 can be calculated by the formulas
 \begin{multline}\label{int-mult-seq-extr-sp-char2}
(h(e^{i\lambda}))^\top=(A_N(e^{i\lambda}))^\top f(\lambda)(f(\lambda)+g(\lambda))^{-1}-\\
-\left(\sum\limits_{j=0}^{N}(\bold{B}_N^{-1}\bold{D}_N\vec{\bold{a}}_N)(j)e^{ij\lambda}\right)^\top(f(\lambda)+g(\lambda))^{-1},
\end{multline}
The mean-square error is of the form
 \[ \Delta (f,g)= \frac{1}{2 \pi} \int_{- \pi}^{ \pi}(A_{N} (e^{i \lambda} )^{ \top} g( \lambda )+C_{N} (e^{i \lambda} )^{ \top} )(f( \lambda )+g( \lambda ))^{-1}  \times \]
 \[ \times f( \lambda )(f( \lambda )+g( \lambda ))^{-1}(A_{N} (e^{i \lambda} )^{ \top} g( \lambda )+C_{N} (e^{i \lambda} )^{ \top} )^{*}d \lambda + \]
 \[+ \frac{1}{2 \pi} \int_{- \pi}^{ \pi}(A_{N} (e^{i \lambda} )^{ \top} f( \lambda )-C_{N} (e^{i \lambda} )^{ \top} )(f( \lambda )+g( \lambda ))^{-1} g( \lambda ) \times \]
 \[ \times(f( \lambda )+g( \lambda ))^{-1} (A_{N} (e^{i \lambda} )^{ \top} f( \lambda )-C_{N} (e^{i \lambda} )^{ \top} )^{*}d \lambda = \]
 \[= \frac{1}{2 \pi} \int_{- \pi}^{ \pi}A_{N} (e^{i \lambda} )^{ \top} f( \lambda )(f( \lambda )+g( \lambda ))^{-1} g( \lambda ) \overline{A_{N} (e^{i \lambda} )}d \lambda + \]
 \[+ \frac{1}{2 \pi} \int_{- \pi}^{ \pi}C_{N} (e^{i \lambda} )^{ \top} (f( \lambda )+g( \lambda ))^{-1} \overline{C_{N} (e^{i \lambda} )}d \lambda = \]
  \begin{equation} \label{int-mult-seq-extr-err2}
  = \left \langle \vec{\bold{a}}_{N} ,\bold{R}_{N} \vec{\bold{a}}_{N} \right \rangle + \left \langle \vec{\bold{c}}_{N} ,\bold{B}_{N} \vec{\bold{c}}_{N} \right \rangle.
 \end{equation}

\begin{thm} \label{theo1.int-stat-noise}
Let $ \vec{ \xi}(j)= \left \{ \xi_{k} (j) \right \}_{k=1}^{T} $ and
$ \vec{ \eta}(j)= \left \{ \eta_{k} (j) \right \}_{k=1}^{T} $ be
uncorrelated vector-valued stationary stochastic sequences with the spectral density matrices
$f( \lambda)= \left\{f_{ij}(\lambda)\right \}_{i,j=1}^{T} $ and
$g( \lambda)=\left\{g_{ij}(\lambda)\right\}_{i,j=1}^{T} $ which satisfy the minimality condition (\ref{int-minimality-stat}).
The spectral characteristic $h(f,g)$ and the value of the mean-square error $ \Delta (f,g)$ of the optimal linear estimate
of the functional $A_N \vec{ \xi}$ from unknown values of stationary stochastic sequence $ \vec{ \xi}(j)$ based on observations of the sequence
$\vec{\xi}(j)+\vec{\eta}(j)$ at points  $j\in\mathbb Z\setminus\{0,1,\dots,N\}$
 are calculated by the formulas (\ref{int-mult-seq-extr-sp-char2}), (\ref{int-mult-seq-extr-err2}).
 \end{thm}

\begin{Corollary}
Let $ \vec{ \xi}(j)= \left \{ \xi_{k} (j) \right \}_{k=1}^{T} $ be a
vector-valued stationary stochastic sequences with spectral density
matrix $f( \lambda)= \left \{f_{ij} ( \lambda ) \right
\}_{i,j=1}^{T}$ which satisfies the minimality condition
 \[ \int_{- \pi}^{ \pi} \,{\rm{Tr}} \, \left[(f( \lambda))^{-1} \right] \,d \lambda< \infty . \]
\noindent The spectral characteristic $h(f)$ and the mean-square error
$ \Delta (f)$ of the optimal linear estimate
of the functional $A_{N} \vec{ \xi}$ from unknown values of stationary stochastic sequence $ \vec{ \xi}(j)$ based on observations of the sequence
$ \vec{ \xi}(j)$ for $j \in Z \backslash \left \{0, \ldots ,N \right \}$
are calculated by the formulas
 \begin{equation} \label{GrindEQ__3_198_}
h(f)^{ \top} =A_{N} (e^{i \lambda} )^{ \top} \, -C_{N} (e^{i \lambda} )^{ \top} [f( \lambda )]^{-1} ,
 \end{equation}
 \begin{equation} \label{GrindEQ__3_199_}
 \Delta (f)= \frac{1}{2 \pi} \int_{- \pi}^{ \pi} \,  C_{N} (e^{i \lambda} )^{ \top} [f( \lambda )]^{-1} \overline{C_{N} (e^{i \lambda} )}d \lambda
 =\left \langle \bold{B}_{N}^{-1} \vec{\bold{a}}_{N} \, , \, \vec{\bold{a}}_{N} \right \rangle ,
 \end{equation}
\noindent where
 \[C_{N} (e^{i \lambda} )= \sum_{j=0}^{N} \vec{c}(j)e^{ij \lambda},\quad \vec{\bold{c}}_{N} =\bold{B}_{N}^{-1} \vec{\bold{a}}_{N} , \]
\noindent $\bold{B}_{N} $ is a matrix constructed from the block-matrices of
dimension $T \times T$:
 \[B_{N} (k, j)=B_{N} (k-j)= \frac{1}{2 \pi} \int_{- \pi}^{ \pi} \left[(f( \lambda ))^{-1} \right]^{ \top} e^{i(j-k) \lambda} d \lambda,\]
 \[ k , j=0, 1, \ldots, N. \]
 \end{Corollary}

 \begin{exm} \label{prikl.1.10}
Consider the problem of estimation of the functional
 \[A{ \kern 1pt}_{1} \vec{ \zeta}=( \alpha \, , \, \beta ) \left( \begin{array}{c} { \zeta_{1} (0)} \\{ \zeta_{2} (0)} \end{array} \right)+( \gamma \, , \, \delta ) \, \left( \begin{array}{c} { \zeta_{1} (1)} \\{ \zeta_{2} (1)} \end{array} \right)= \]
 \[= \alpha \, \zeta_{1} (0)+ \beta \, \zeta_{2} (0)+ \gamma \, \zeta_{1} (1)+ \delta \, \zeta_{1} (1) \]
from the unknown values of two-dimensional stationary sequence
$ \vec{ \zeta}(n)= \left \{ \zeta_{k} (n) \right \}_{k=1}^{2}$ based on observations
$ \vec{ \zeta}(j)$, $j \in Z \backslash \left \{0 \, , \, 1 \right \}$,
where $ \zeta_{1} (n)= \xi (n)$ is a stationary stochastic sequence with the  spectral density $f( \lambda )$, and $ \zeta_{2} (n)= \xi (n)+ \eta (n)$, where $ \eta (n)$ is an uncorrelated with $ \xi (n)$ stationary stochastic sequence with the  spectral density $g( \lambda )$. The matrix of spectral densities is of the form
 \[F( \lambda )= \left( \begin{array}{cc} {f( \lambda )} & {f( \lambda )} \\{f( \lambda )} & {f( \lambda )+g( \lambda )} \end{array} \right). \]
\noindent Its determinant equals
 \[D= \left|F( \lambda ) \right|=f( \lambda ) \, g( \lambda ), \]
\noindent and the inverse matrix is as follows
 \[F( \lambda )^{-1} = \left( \begin{array}{cc} { \frac{f( \lambda )+g( \lambda )}{f( \lambda ) \, g( \lambda )}} & { \frac{-1}{g( \lambda )}} \\{ \frac{-1}{g( \lambda )}} & { \frac{1}{g( \lambda )}} \end{array} \right). \]
\noindent Let
 \[f( \lambda )= \frac{P_{1}}{2 \pi \, \left|1-b \, e^{i \lambda} \right|^{2}}, \quad g( \lambda )= \frac{P_{2}}{2 \pi} , \quad P_{1} , P_{2} , b \in R, |b|<1. \]
\noindent Then
 \[ \frac{f( \lambda )+g( \lambda )}{f( \lambda ) \, g( \lambda )} = \frac{2 \pi}{P_{1} P_{2}} \left(P_{1} +P_{2} +P_{2} b^{2} -P_{2} b e^{-i \lambda} -P_{2} b e^{i \lambda} \right),
 \]
 \[ \frac{1}{g( \lambda )} = \frac{2 \pi}{P_{2}} . \]

\noindent We have the matrix $B_{1} $:
 \[B_{1} = \left( \begin{array}{cccc} { \frac{2 \pi}{P_{1} \, P_{2}} (P_{1} +P_{2} +P_{2} b^{2} )} & {- \frac{2 \pi}{P_{2}}} & {- \frac{2 \pi}{P_{1}} b} & {0} \\{- \frac{2 \pi}{P_{2}}} & { \frac{2 \pi}{P_{2}}} & {0} & {0} \\{- \frac{2 \pi}{P_{1}} b} & {0} & { \frac{2 \pi}{P_{1} \, P_{2}} (P_{1} +P_{2} +P_{2} b^{2} )} & {- \frac{2 \pi}{P_{2}}} \\{0} & {0} & {- \frac{2 \pi}{P_{2}}} & { \frac{2 \pi}{P_{2}}} \end{array} \right) \]
\noindent Its determinant equals
 \[D= \left( \frac{2 \pi}{P_{2}} \right)^{4} \cdot \frac{P_{2}^{2}}{P_{1}^{2}} (1+b^{2} +b^{4} ). \]
\noindent The inverse matrix $B_{1}^{-1}$ is equal to
 \[ \frac{P_{2}}{D \cdot P_{1}} \left( \begin{array}{cccc} {1+b^{2}} & {1+b^{2}} & {b} & {b} \\{1+b^{2}} & { \frac{P_{1} +P_{2}}{P_{1}} (1+b^{2} )+ \frac{P_{2}}{P_{1}} b^{4}} & {b} & {b} \\{b} & {b} & {1+b^{2}} & {1+b^{2}} \\{b} & {b} & {1+b^{2}} & { \frac{P_{1} +P_{2}}{P_{1}} (1+b^{2} )+ \frac{P_{2}}{P_{1}} b^{4}} \end{array} \right) \]
\noindent The vector $ \vec{c}_{1}$ is as follows
 \[ \vec{c}_{1} =B_{1}^{-1} \vec{a}_{1} =B_{1}^{-1} \left( \begin{array}{c} { \alpha} \\{ \beta} \\{ \gamma} \\{ \delta} \end{array} \right)=
 \frac{1}{2 \pi (1+b^{2} +b^{4} )} \times
 \]
 \[\times\left( \begin{array}{c} {P_{1} \left[( \alpha + \beta )(1+b^{2} )+b( \gamma + \delta ) \right]} \\{P_{1} \left[( \alpha + \beta )(1+b^{2} )+b( \gamma + \delta ) \right]+P_{2} { \kern 1pt} \beta { \kern 1pt} (1+b^{2} +b^{4} )} \\{P_{1} \left[(1+b^{2} )( \gamma + \delta )+b( \alpha + \beta ) \right]} \\{P_{1} \left[b( \alpha + \beta )+( \gamma + \delta )(1+b^{2} ) \right]+P_{2} { \kern 1pt} \delta { \kern 1pt} (1+b^{2} +b^{4} )} \end{array} \right) \]
\noindent The spectral characteristic of the optimal estimate of the random variable $A_{1} \vec{ \zeta}$ can be calculated by the formula
 \[h(F)^{ \top} =( \alpha + \gamma \, e^{i \lambda} \, , \, \beta + \delta \, e^{i \lambda} )-C_{1} (e^{i \lambda} )^{ \top} F( \lambda )^{-1} =(h_{1} \,, h_{2}), \]
\noindent where
 \[h_{1}=A \left[( \alpha + \beta )(1+b^{2} )+b( \gamma + \delta ) \right] e^{-i \lambda}+ \]
 \[+A \left[( \gamma + \delta )(1+b^{2} )+b( \alpha + \beta ) \right] e^{2i \lambda}, \, h_{2}=0, \,A= \frac{b}{1+b^{2} +b^{4}} . \]
\noindent Thus the optimal estimate of the random variable $A_{1} \vec{ \zeta}$ is of the form
 \[ \hat{A}_{1} \vec{ \zeta}=A \left[( \alpha + \beta )(1+b^{2} )+b( \gamma + \delta ) \right] \, \zeta_{1} (-1)+ \]
 \[+ A \left[( \gamma + \delta )(1+b^{2} )+b( \alpha + \beta ) \right] \, \zeta_{1} (2). \]
\noindent The value of the mean-square error equals
\[\Delta (F)=\left\langle \vec{c}_{1} \, ,\, \vec{a}_{1} \right\rangle =\]
\[=\frac{AP_{1} }{2\pi b} \left[(\alpha +\beta )^{2} (1+b^{2} )+(\gamma +\delta )^{2} (1+b^{2} )+2b(\gamma +\delta )(\alpha +\beta )\right]+\]
\[+\frac{P_{2} }{2\pi } \left[\beta ^{2} +\delta ^{2} \right].\]
\end{exm}

\medskip

 \begin{exm} \label{prikl.1.11}
Let all conditions of the previous example are satisfied and let
 \[f( \lambda )= \frac{P_{1}}{2 \pi \, \left|1-b_{1} e^{i \lambda} \right|^{2}} , \quad g( \lambda )= \frac{P_{2}}{2 \pi \, \left|1-b_{2} e^{i \lambda} \right|^{2}},\]
 $
 P_{1}, \, P_{2}, \, b_{1},b_{2} \in R,  |b_1|<1, |b_2|<1. $

\noindent Then we will have
 \[B_{1}^{-1} = \left( \begin{array}{cccc} { \frac{P_{1}}{2 \pi} A} & { \frac{P_{1}}{2 \pi} A} & { \frac{P_{1}}{2 \pi} C} & { \frac{P_{1}}{2 \pi} C} \\{ \frac{P_{1}}{2 \pi} A} & { \frac{P_{1}}{2 \pi} A+ \frac{P_{2}}{2 \pi} B} & { \frac{P_{1}}{2 \pi} C} & { \frac{P_{1}}{2 \pi} C+ \frac{P_{2}}{2 \pi} D} \\{ \frac{P_{1}}{2 \pi} C} & { \frac{P_{1}}{2 \pi} C} & { \frac{P_{1}}{2 \pi} A} & { \frac{P_{1}}{2 \pi} A} \\{ \frac{P_{1}}{2 \pi} C} & { \frac{P_{1}}{2 \pi} C+ \frac{P_{2}}{2 \pi} D} & { \frac{P_{1}}{2 \pi} A} & { \frac{P_{1}}{2 \pi} A+ \frac{P_{2}}{2 \pi} B} \end{array} \right), \]
\noindent where
 \[A= \frac{1+b_{1}^{2}}{1+b_{1}^{2} +b_{1}^{4}},\quad B= \frac{1+b_{2}^{2}}{1+b_{2}^{2} +b_{2}^{4}},
\]
\[ C= \frac{b_{1}}{1+b_{1}^{2} +b_{1}^{4}}, \quad D= \frac{b_{2}}{1+b_{2}^{2} +b_{2}^{4}} . \]
The spectral characteristic of the optimal estimate of the random variable $A_{1} \vec{ \zeta}$ can be calculated by the formula
 \[h(F)^{ \top} =( \alpha + \gamma \, e^{i \lambda} \, , \, \beta + \delta \, e^{i \lambda} )-C_{1} (e^{i \lambda} )^{ \top} F( \lambda )^{-1} =(h_{{ \kern 1pt} 1} (e^{i \lambda} ),h_{2} (e^{i \lambda} )), \]
where
\begin{multline*}
h_{1} (e^{i \lambda} )=\\
=(C \left[( \alpha + \beta )(1+b_{1}^{2} )+b_{1} ( \gamma + \delta )
\right] \, -D \left[ \beta (1+b_{2}^{2} )+b_{2} \delta \right] \, ){
\kern 1pt} { \kern 1pt} e^{-i \lambda} +
\end{multline*}
 \[+(C \left[( \gamma + \delta )(1+b_{1}^{2} )+b_{1} ( \alpha + \beta ) \right]-D \left[ \delta (1+b_{2}^{2} )+b_{2} \beta \right] \, ) \, e^{2i \lambda} ; \]
 \[h_{2} (e^{i \lambda} )=D \left[ \beta (1+b_{2}^{2} )+b_{2} \delta \right]{ \kern 1pt} { \kern 1pt} e^{-i \lambda} +D \left[ \delta (1+b_{2}^{2} )+b_{2} \beta \right] \, ){ \kern 1pt} { \kern 1pt} e^{2i \lambda} . \]
The optimal estimate of the random variable $A_{1} \vec{ \zeta}$ is of the form
\begin{multline*}
\hat{A}_{1} \vec{ \zeta}=\\=(C \left[( \alpha + \beta )(1+b_{1}^{2}
)+b_{1} ( \gamma + \delta ) \right] \, -D \left[ \beta (1+b_{2}^{2}
)+b_{2} \delta \right] \, ){ \kern 1pt} { \kern 1pt} \zeta_{1} (-1)+
\end{multline*}
 \[+(C \left[( \gamma + \delta )(1+b_{1}^{2} )+b_{1} ( \alpha + \beta ) \right]-D \left[ \delta (1+b_{2}^{2} )+b_{2} \beta \right] \, ) \, \zeta_{1} (2)+ \]
 \[+D \left[ \beta (1+b_{2}^{2} )+b_{2} \delta \right]{ \kern 1pt} { \kern 1pt} \zeta_{2} (-1)+D \left[ \delta (1+b_{2}^{2} )+b_{2} \beta \right] \, \zeta_{2} (2). \]
The value of the mean-square error is
 \[ \Delta (F)= \left \langle \vec{c}_{1} \, , \, \vec{a}_{1} \right \rangle = \frac{P_{1}}{2 \pi (1+b_{1}^{2} +b_{1}^{4} )} \times \]
 \[ \times \left[( \alpha + \beta )^{2} (1+b_{1}^{2} )+( \gamma + \delta )^{2} (1+b_{1}^{2} )+2b_{1} ( \gamma + \delta )( \alpha + \beta ) \right]+ \]
 \[+ \frac{P_{2}}{2 \pi (1+b_{2}^{2} +b_{2}^{4} )} \left[( \beta ^{2} + \delta ^{2} )(1+b_{2}^{2} )+2b_{2} \delta \beta \right]. \]
 \end{exm}

 \section{Interpolation problem. Minimax approach}

  The value of the error
 $$\Delta\left(h(f_0,g_0);f,g\right):= \left\| {A}_N \vec{\xi}- \hat{A}_N  \vec{\xi}\,\, \right\|_\alpha^\alpha $$
 and
the spectral characteristic $h(f,g):={h}(\theta)$ of the optimal estimate
$\hat{A}_N \vec{\xi}$ of the functional $A_N \vec{\xi} =\sum_{j = 0}^{N} (\vec{a}(j))^{\top} \vec{\xi}(j)$
that depends on the unknown values $\vec{\xi}(j),j=0,1,\dots,N,$  of the sequence $\vec{\xi}(j)$, from observations of the sequence
 $\vec{\xi}(j)+\vec{\eta}(j)$  at points  $j\in\mathbb Z\setminus\{0,1,\dots,N\}$
can be calculated by the proposed formulas
 only in the case where we know the spectral densities $f(\theta)$ and $g(\theta)$ of the
 mutually independent harmonizable symmetric $\alpha$-stable $HS{\alpha}S$, stochastic sequences $\vec{\xi}(j)=\left \{ \xi_ {k} (j) \right \}_{k = 1} ^ {T},\,j\in \mathbb{Z}$  and
 $\vec{\eta}(j)=\left \{ \eta_ {k} (j) \right \}_{k = 1} ^ {T},\,j\in \mathbb{Z}$.

 However, in practice we can't exactly evaluate the spectral densities of stochastic sequences, but, instead,  we often can have a set $D=D_f\times D_g$ of admissible spectral densities. In this case we can apply the minimax-robust method of estimation to the interpolation problem. This method let us find an estimate that minimizes
the maximum of the errors for all spectral densities from the
given set $D=D_f\times D_g$ of admissible spectral
densities simultaneously (see books by Moklyachuk (2008a), Moklyachuk and Masyutka (2012), Moklyachuk and Golichenko (2016)).

 \begin{definition}
 For a given class of spectral densities $D=D_f\times D_g$ the spectral densities $f_0(\theta)\in D_f$, $g_0(\theta)\in D_g$ are called the least favorable in $D=D_f\times D_g$ for the optimal linear estimation $\hat{A}_N \vec{\xi}$ of the functional ${A}_N \vec{\xi}$, if the following relation holds true
 $$\Delta\left(f_0,g_0\right)=\Delta\left(h\left(f_0,g_0\right);f_0,g_0\right)=\max\limits_{(f,g)\in D_f\times D_g}\Delta\left(h\left(f,g\right);f,g\right).$$
 \end{definition}

 \begin{definition}
 For a given class of spectral densities $D=D_f\times D_g$ the spectral characteristic $h^0=h(f_0)$
 of the optimal estimate $\hat{A}_N \vec{\xi}$ of the functional ${A}_N \vec{\xi}$ is called
 minimax (robust)
 for the optimal linear estimation ${A}_N \vec{\xi}$, if the following relations hold true
$$h^0(\theta)\in H_D= \bigcap\limits_{(f,g)\in D_f\times D_g} L^{\alpha}(f+g),$$
$$\min\limits_{h\in H_D}\max\limits_{(f,g)\in D}\Delta\left(h;f,g\right)=\max\limits_{(f,g)\in D}\Delta\left(h^0;f,g\right).$$
 \end{definition}

 The least favorable spectral densities $f_0(\theta)$, $g_0(\theta)$ and the minimax spectral characteristic $h^0=h(f_0,g_0)$ form a saddle point of the function $\Delta \left(h;f,g\right)$ on the set $H_D \times D$. The saddle point inequalities
$$\Delta\left(h;f_0,g_0\right)\geq\Delta\left(h^0;f_0,g_0\right)\geq \Delta\left(h^0;f,g\right) $$ $$ \forall h \in H_D, \forall f \in D_f, \forall g \in D_g$$
holds true if $h^0=h(f_0,g_0)$ and $h(f_0,g_0)\in H_D,$ where $(f_0,g_0)$ is a solution to the conditional extremum problem
\begin{equation} \label{condextr}
\max\limits_{(f,g)\in D_f\times D_g}\Delta\left(h(f_0,g_0);f,g\right)=\Delta\left(h(f_0,g_0);f_0,g_0\right),
\end{equation}
where
\begin{multline}\label{delta41}
\Delta\left(h(f_0,g_0);f,g\right)=\left\|{A}_N \vec{\xi}- \hat{A}_N \vec{\xi} \right\|_\alpha^\alpha =
\\
=\int_{-\pi}^{\pi}
\left( A_N(e^{i\theta}) - {h}(\theta) \right)^{\top}
f(\theta)
\left( A_N(e^{i\theta}) - {h}(\theta) \right)^{<{\alpha - 1}>}
d \theta+
\\
 +
 \int_{-\pi}^{\pi}
\left( {h}(\theta) \right)^{\top}
g(\theta)
\left(  {h}(\theta) \right)^{<{\alpha - 1}>}
d \theta
 \end{multline}

The constrained optimization problem (\ref{condextr}) is equivalent to the unconditional extremum problem
\begin{equation} \label{8}
\Delta_D(f,g)=-\Delta(h(f_0,g_0);f,g)+\delta(f,g\left|D_f\times D_g\right.)\rightarrow \inf,
\end{equation}
where $\delta(f,g\left|D_f\times D_g\right.)$ is the indicator function of the set $D=D_f\times D_g$.
 Solution $(f_0,g_0)$ to the problem (\ref{8}) is characterized by the condition $0 \in \partial\Delta_D(f_0,g_0),$
 where $\partial\Delta_D(f_0)$ is the subdifferential of the convex functional $\Delta_D(f,g)$ at point $(f_0,g_0)$.
This condition makes it possible to find the least favorable spectral densities in some special classes of spectral densities
$D=D_f\times D_g$ (Ioffe and Tihomirov, 1979; Pshenychnyj, 1971; Rockafellar, 1997; Moklyachuk, 2008b).

Note, that the form of the functional $\Delta(h(f_0,g_0);f,g)$ is convenient for application the Lagrange method of indefinite multipliers for
finding solution to the problem (\ref{condextr}).
Making use the method of Lagrange multipliers and the form of
subdifferentials of the indicator functions
we describe relations that determine least favourable spectral densities in some special classes
of spectral densities

Summing up the derived formulas and the introduced definitions we come to conclusion that the following lemmas hold true

\begin{lemma} \label{lem41}
Let $\vec{\xi}(j)=\left \{ \xi_ {k} (j) \right \}_{k = 1} ^ {T},\,j\in \mathbb{Z}$  and
 $\vec{\eta}(j)=\left \{ \eta_ {k} (j) \right \}_{k = 1} ^ {T},\,j\in \mathbb{Z}$ be mutually independent harmonizable symmetric $\alpha$-stable random sequences which have
absolutely continuous spectral measures and the spectral densities $f_0(\theta)$ and $g_0(\theta)$ satisfying the minimality condition (\ref{minimalityF+G}).
Let the spectral densities $(f_0,g_0)\in D_f\times D_g$ gives a solution to the constrained optimization problem  (\ref{condextr}).
The spectral densities $(f_0,g_0)$ are the least favorable spectral densities in $D_f\times D_g$ and
$h^0=h(f_0,g_0)$ is the minimax spectral characteristic of the optimal linear estimate
$\hat{A}_N \vec{\xi}$ of the functional $A_N \vec{\xi}$
that depends on the unknown values $\vec{\xi}(j),j=0,1,\dots,N,$  of the sequence $\vec{\xi}(j)$, from observations of the sequence
$\vec{\xi}(j)+\vec{\eta}(j)$  at points  $j\in\mathbb Z\setminus\{0,1,\dots,N\}$,
if $h^0=h(f_0,g_0)\in H_D$.
\end{lemma}

\begin{lemma} \label{lem42}
Let $\vec{\xi}(j)=\left \{ \xi_ {k} (j) \right \}_{k = 1} ^ {T},\,j\in \mathbb{Z}$  be a harmonizable symmetric $\alpha$-stable random sequence which has absolutely continuous spectral measure
and the spectral density $f_0(\theta)$ satisfying the minimality condition
(\ref{minimality}).
Let the spectral density $f_0\in D_f$ gives a solution to the constrained optimization problem
\begin{equation} \label{condextr42}
\max\limits_{f\in D_f}\Delta\left(h(f_0);f\right)=\Delta\left(h(f_0);f_0\right),
\end{equation}
where
   \begin{multline}\label{condextr43}
\Delta\left(h(f_0);f\right)=\left\| {A}_N \vec{\xi}- \hat{A}_N \vec{\xi}\,\, \right\|_\alpha^\alpha =
\\
=
\int_{-\pi}^{\pi}
\left[\left(f^{-1}(\theta)C_N(e^{i\theta}) \right)^{<\frac{1}{\alpha - 1}>}\right]^{\top}
f(\theta)
\left(f^{-1}(\theta)C_N(e^{i\theta}) \right)^{<\frac{\alpha-1}{\alpha - 1}>}
d \theta.
\end{multline}
The spectral density $f_0$ is the least favorable spectral density in $D_f$ and
$h^0=h(f_0)$ is the minimax spectral characteristic of the optimal linear estimate $\hat{A}_N \vec{\xi}$ of the functional $A_N \vec{\xi}$
that depends on the unknown values $\vec{\xi}(j),j=0,1,\dots,N,$  of the sequence $\vec{\xi}(j)$, from observations of the sequence
$\vec{\xi}(j)$  at points  $j\in\mathbb Z\setminus\{0,1,\dots,N\}$,
if $h^0=h(f_0)\in H_D$.
\end{lemma}

     \subsection{Least favorable spectral densities in the class  $D_f^0$}

Consider the problem of optimal linear estimation of the functional
 $A_N \xi= \sum_{j = 0}^{N} a(j) \xi(j)$ that depends on the
unknown values $\xi(j), j = 0, 1, \ldots,N$, of a random sequence $\{\xi(k),k\in\mathbb Z\},$ from observations of the sequence $\{\xi(j),j\in\mathbb Z\}$ at points  $j\in\mathbb Z\setminus\{0,1,\dots,N\}$,
where $\{\xi(k),k\in\mathbb Z\}$ is a harmonizable symmetric $\alpha$-stable random sequence which have the
spectral density $f_0(\theta)$ satisfying the minimality condition (\ref{minimality}) from
 the class of admissible spectral densities $D_0$, where
$$D_f^0= \left\{f: \int_{-\pi}^{\pi} f(\theta) d \theta = P \right\}.$$
Applying the Lagrange multipliers method  to find a solution to the constrained optimization problem (\ref{condextr42}) we derive that
the least favorable density $f_0\in D_f^0$ satisfy the equation
       $$\int_{-\pi}^{\pi}\left( \left|\left(\sum_{j = 0}^{N} c(j) e^{-i j \theta} \right)^{<\frac{1}{\alpha - 1}>} (f_0(\theta))^{\frac{-1}{\alpha - 1}} \right|^{\alpha} - \lambda \right)f(\theta) d \theta = 0.$$
It follows from the Lagrange lemma that
    $$\left|\left(\sum_{j = 0}^{N} c(j) e^{-i j \theta} \right)^{<\frac{1}{\alpha - 1}>} (f_0(\theta))^{\frac{-1}{\alpha - 1}} \right|^{\alpha} =\lambda.$$
The least  favorable density $f_0\in D_f^0$  is of the form
\begin{equation}\label{least1}
f_0(\theta) = C\left|\sum_{j = 0}^{N} c(j) e^{-i j \theta} \right|.
\end{equation}
Thus, the following statement holds true.

\begin{thm}
The spectral density (\ref{least1}), which satisfy the minimality condition (\ref{minimality}),
equation (\ref{sp_eq2}) and condition $f_0\in D_0$, that is $\int_{-\pi}^{\pi} f_0(\theta) d \theta = P$,
is the least  favorable density $f_0\in D_f^0$
for the optimal linear estimation of the functional
 $A_N \xi= \sum_{j = 0}^{N} a(j) \xi(j)$.
 The minimax-robust spectral characteristic $h(f_0)$ of the optimal estimate of the
functional $A_N \xi$ is determined by the
 (\ref{spectralcharacteristic})  with $f(\theta)=f_0(\theta)$.
\end{thm}

     \subsection{Least favorable spectral densities in the class $D_f^{\beta}$}

Consider the problem of optimal linear estimation of the functional
 $A_N \xi= \sum_{j = 0}^{N} a(j) \xi(j)$ that depends on the
unknown values $\xi(j), j = 0, 1, \ldots,N$, of a random sequence $\{\xi(k),k\in\mathbb Z\},$ from observations of the sequence $\{\xi(j),j\in\mathbb Z\}$ at points  $j\in\mathbb Z\setminus\{0,1,\dots,N\}$,
where $\{\xi(k),k\in\mathbb Z\}$ is a harmonizable symmetric $\alpha$-stable random sequence which have the
spectral density $f_0(\theta)$ satisfying the minimality condition (\ref{minimality}) from
 the class of admissible spectral densities
$$D_f^{\beta} = \left\{f: \int_{-\pi}^{\pi} \left(f(\theta)\right)^{\beta} d \theta = P \right\}, \quad \beta \neq \frac{-1}{\alpha - 1}, \beta \neq 1.$$

Applying the Lagrange multipliers method  to find a solution to the constrained optimization problem (\ref{condextr42}) we derive that
the least favorable density $f_0\in D_f^{\beta}$ satisfy the equation

$$\left|\left(\sum_{j = 0}^{N} c_j e^{-i j \theta} \right)^{<\frac{1}{\alpha - 1}>} (f_0(\theta))^{\frac{-1}{\alpha - 1}} \right|^{\alpha} = \lambda \left(f_0(\theta)\right)^{\beta - 1},$$
which can be written in the form
$$\left|\sum_{j = 0}^{N} c_j e^{-i j \theta} \right|^{\frac{\alpha}{\alpha - 1}} (f_0(\theta))^{\frac{-\alpha}{\alpha - 1}} = \lambda \left(f_0(\theta)\right)^{\beta - 1}.$$
It follows from this relation that
the least favorable density $f_0\in D_f^{\beta}$ is of the form
\begin{equation}\label{least2}
f_0(\theta) = C\left|\sum_{j = 0}^{N} c_j e^{-i j \theta} \right|^\frac{-\alpha}{-\alpha - (\alpha-1)(\beta - 1)}.
\end{equation}

Thus, the following statement holds true.

\begin{thm}
The spectral density (\ref{least2}), which satisfy the minimality condition (\ref{minimality}),
equation (\ref{sp_eq2}) and condition $f_0\in D_f^{\beta}$, that is $\int_{-\pi}^{\pi} \left(f_0(\theta)\right)^{\beta} d \theta = P$,
is the least  favorable density $f_0\in D_f^{\beta}$
for the optimal linear estimation of the functional
 $A_N \xi= \sum_{j = 0}^{N} a(j) \xi(j)$.
 The minimax-robust spectral characteristic $h(f_0)$ of the optimal estimate of the
functional $A_N \xi$ is determined by the
 (\ref{spectralcharacteristic})  with $f(\theta)=f_0(\theta)$.
\end{thm}

     \subsection{Least favorable spectral densities in the class $D_{M}^{-}$}

Consider the problem of optimal linear estimation of the functional
 $A_N \xi= \sum_{j = 0}^{N} a(j) \xi(j)$ that depends on the
unknown values $\xi(j), j = 0, 1, \ldots,N$, of a random sequence $\{\xi(k),k\in\mathbb Z\},$ from observations of the sequence $\{\xi(j),j\in\mathbb Z\}$ at points  $j\in\mathbb Z\setminus\{0,1,\dots,N\}$,
where $\{\xi(k),k\in\mathbb Z\}$ is a harmonizable symmetric $\alpha$-stable random sequence which have the
spectral density $f_0(\theta)$ satisfying the minimality condition (\ref{minimality}) from
 the class of admissible spectral densities
$$D_{M}^{-} = \left\{f: \int_{-\pi}^{\pi} f^{-1}(\theta) cos(m \theta) d \theta = r_m, m = 0, \ldots, M  \right\},$$
where
 $r_m,$ $m = 0, \ldots, M$ is a strictly positive sequence of real numbers.
 Under this condition the moment problem has solutions and the set $D_{M}^{-}$ contains an infinite number of densities (Krein and Nudelman, 1977).

Applying the Lagrange multipliers method  to find a solution to the constrained optimization problem (\ref{condextr42}) we derive that
the least favorable density $f_0\in D_f^{\beta}$ satisfy the equation
$$\left|\left(\sum_{j = 0}^{N} c_j e^{-i j \theta} \right)^{<\frac{1}{\alpha - 1}>} (f_0(\theta))^{\frac{-1}{\alpha - 1}} \right|^{\alpha} =
\left(\sum_{m = 0}^{M} \lambda_{m}  \cos(m \theta)\right)(f_0(\theta))^{-2},$$
where $\lambda_j,j=0, \ldots, M $ are Lagrange multipliers.
The function
$$\left(\sum_{m = 0}^{M} \lambda_{m}  \cos(m \theta)\right)$$
can be represented in the form (Hannan, 1970)
$$
\left(\sum_{m = 0}^{M} \lambda_{m}  \cos(m \theta)\right)= \left|\sum_{m = 0}^{M} p_{m} e^{im \theta}\right|^2.
$$
It follows from the obtained relations that
the least favorable density $f_0\in D_{M}^{-}$ is of the form
($1<\alpha < 2$)
\begin{equation}\label{least3}
f_0(\theta) = \left|\sum_{j = 0}^{N} c_j e^{-i j \theta} \right|^{\frac{\alpha}{2 - \alpha}}
{\left|\sum_{m = 0}^{M} p_{m} e^{im \theta}\right|^{2\frac{1 - \alpha}{2 - \alpha}}}.
\end{equation}

Thus, the following statement holds true.

\begin{thm}
The spectral density (\ref{least3}), which satisfy the minimality condition (\ref{minimality}),
equation (\ref{sp_eq2}) and condition $f_0\in D_{M}^{-}$, that is $\int_{-\pi}^{\pi} f_0^{-1}(\theta) cos(m \theta) d \theta = r_m, m = 0, \ldots, M$,
is the least  favorable density $f_0\in D_{M}^{-}$
for the optimal linear estimation of the functional
 $A_N \xi= \sum_{j = 0}^{N} a(j) \xi(j)$.
 The minimax-robust spectral characteristic $h(f_0)$ of the optimal estimate of the
functional $A_N \xi$ is determined by the
 (\ref{spectralcharacteristic})  with $f(\theta)=f_0(\theta)$.
\end{thm}

\subsection{Least favorable spectral densities in the class $D_f^{-1}$}

Consider the problem of optimal linear estimation of the functional
 $A_N \xi= \sum_{j = 0}^{N} a(j) \xi(j)$ that depends on the
unknown values $\xi(j), j = 0, 1, \ldots,N$, of a random sequence $\{\xi(k),k\in\mathbb Z\},$ from observations of the sequence $\{\xi(j),j\in\mathbb Z\}$ at points  $j\in\mathbb Z\setminus\{0,1,\dots,N\}$,
where $\{\xi(k),k\in\mathbb Z\}$ is a harmonizable symmetric $\alpha$-stable random sequence which have the
spectral density $f_0(\theta)$ satisfying the minimality condition (\ref{minimality}) from
 the class of admissible spectral densities
 $D_f^{\beta}$, where $\beta=-1,$ and the sequence $a(j), j=0,1,\dots,N,$ that determines the functional $A\xi$, is strictly positive.

By using the method of Lagrange multipliers we get that the Fourier coefficients of the function $f_0^{-1}$ satisfy the equation
\begin{equation} \label{99}
\left|\sum\limits_{j=0}^{N}c(j)e^{ij\theta} \right|^2=p_0^2,
\end{equation}
where $c(j), j=0,1,\dots,N,$ are components of the vector $\vec{\bold{c}}_N=\{c(j): j=0,1,\dots,N,\}$ that satisfies the equation $\bold{B}_N^0\vec{\bold{c}}_N=\vec{\bold{a}}_N$,
$\vec{\bold{a}}_N=\{a(j): j=0,1,\dots,N,\}$,
the matrix $\bold{B}_N^0$ is determined by the Fourier coefficients of the function $f_0^{-1}(\theta)$
\begin{equation*}
B_N^0(k,j)=\,\frac{1}{2\pi}\int\limits_{-\pi}^{\pi}f_0^{-1}(\theta)e^{-i(k-j)\theta}d\theta=r_{k-j}^0,
\end{equation*}
\begin {equation*}
k,j=0,1,\dots,N.
\end{equation*}
The Fourier coefficients $r_k=r_{-k}, k=0,1,\dots,N,$ satisfy both equation (\ref{99}) and equation $\bold{B}_N^0\vec{\bold{c}}_N=\vec{\bold{a}}_N$.
These coefficients can be found from the equation
 $\bold{B}_N^0\vec{\bold{p}}_N^0=\vec{\bold{a}}_N$,
 where $\bold{p}_N^0=(p_0,0,\ldots,0, \ldots).$
 The last relation can be presented in the form of the system of equations
\begin {eqnarray*}
&r_k p_0=a(k), \, k=0,1,\dots,N.
\end{eqnarray*}
From the first equation of the system (for $k=0$) we find the unknown value $p_0=a(0) (r_0)^{-1}$.
 It follows from the restriction on the spectral densities from the class $D_0^{-1}$ that the Fourier coefficient
 $$r_0=\int\limits_{-\pi}^{\pi}f_0^{-1}(\theta)d\theta=P_1.$$
 For this reason
 $$r_k=r_{-k}=P_1a(k) (a(0))^{-1}, k=1,\dots,N.$$
We can represent the function $f_0^{-1}(\theta)$ in the form
 $$f_0^{-1}(\theta)=\sum_{k=-N}^N r_k e^{ik\theta}.$$
 Since the sequence $a(j), j=0,1,\dots,N,$ is strictly positive,
 the sequence $r_k, k= 0,1,\ldots,N$, is also strictly positive and the function $f_0^{-1}(\theta)$ is positive,
 so it can be represented in the form (Hannan, 1970; Krein and Nudelman, 1977)
 \[f_0^{-1}(\theta)=\left|\sum\limits_{k=0}^{N}\gamma_k e^{-ik\theta}\right|^2,\,\, \theta \in \left[-\pi, \pi\right].
 \]
Hence, $f_0(\theta)$ is the spectral density of the autoregressive stochastic sequence of the infinite order generated by equation
\begin{equation} \label{regr}
\sum\limits_{k=0}^{N}\gamma_k \xi(n-k)=\epsilon_n,
\end{equation}
where $\epsilon_n$ is a ``white noise'' sequence.

Thus, the following theorem holds true.

\begin{thm}
The least favorable in the class $D_0^{-1}$ spectral density for the optimal linear estimation of the functional
$A\xi$ determined by strictly positive sequence $a(j), j=0,1,\dots,N,$
is the spectral density of the autoregressive sequence (\ref{regr}) whose Fourier coefficients are $r_k=r_{-k}=P_1a(k) (a(0))^{-1}, k=0,1,\dots,N.$
\end{thm}

 \section{Conclusion}

 We propose methods of solution the optimal linear estimation problem for the linear functionals
that depend on the
unknown values  $\vec{\xi}(j),j=0,1,\dots,N,$  of a vector-valued harmonizable symmetric $\alpha$-stable random sequence $\vec{\xi}(j)=\left \{ \xi_ {k} (j) \right \}_{k = 1} ^ {T}$,
from observations of the sequence $\vec{\xi}(j)+\vec{\eta}(j)$ at points  $j\in\mathbb Z\setminus\{0,1,\dots,N\}$
where $\vec{\xi}(j)$ and $\vec{\eta}(j)$ are mutually independent harmonizable symmetric $\alpha$-stable random sequences
which have the spectral densities $f(\theta)$ and $g(\theta)$ satisfying the minimality condition.

The problem is investigated under the condition of spectral certainty as well as under the condition of spectral uncertainty.
 Formulas for calculation the value of the error and spectral characteristic of the optimal linear
estimate of the functional are derived under the condition of spectral certainty where spectral density of the sequence is exactly known.
 In the case where spectral density of
the sequence is not exactly known, but a set of admissible spectral densities is available, relations which determine least favorable densities and the minimax-robust spectral characteristics for different classes of spectral densities are found.


\newpage
 \begin{thebibliography}{99}





\bibitem{Cambanis}
 Cambanis, S., 1983.
 Complex stable variables and processes, Contributions to Statistics: Essays in Honour
of Norman L. Johnson, P. K. Sen, ed., North-Holland, New York, 63-79.


\bibitem{CambanisSoltani1984}
 Cambanis, S., Soltani, R., 1984.
 Prediction of Stable Processes:
Spectral and Moving Average Representations.
Z. Wahrscheinlichkeitstheorie verw. Gebiete. 66, 593-612.




\bibitem{Franke}
Franke, J., 1985. Minimax robust prediction of discrete
time series. Z. Wahrsch. Verw. Gebiete 68, 337-364.

\bibitem{FrankePoor}
Franke, J., Poor, H. V., 1984. Minimax-robust filtering
and finite-length robust predictors. Robust and Nonlinear Time
Series Analysis. Lecture Notes in Statistics, Springer-Verlag
26, 87-126.


\bibitem{Grenander}
Grenander, U., 1957.
A prediction problem in game theory. Ark. Mat. 3, 371-379.

\bibitem{Hannan}
Hannan, E.J., 1970.
Multiple time series. Wiley Series in Probability and Mathematical Statistics. New York etc.: John Wiley \& Sons, Inc. XI, 536.


\bibitem{Hosoya}
Hosoya, Y., 1982. Harmonizable stable processes. Z. Wahrsch. Verw. Gebiete 60, 517-533.

\bibitem{Ioffe}
Ioffe, A. D., Tihomirov, V. M., 1979.
Theory of extremal problems. Amsterdam, New York, Oxford:
North--Holland Publishing Company, 460.

\bibitem{Kassam}
Kassam, S.A., Poor, H.V., 1985.
Robust techniques for signal processing: A survey. Proceedings of the IEEE 73,
433-481.

\bibitem{Kolmogorov}
Kolmogorov, A.N., 1992. Selected works of A. N. Kolmogorov. Vol. II: Probability theory and mathematical
statistics. Ed. by A. N. Shiryayev. Mathematics and Its Applications. Soviet Series. 26. Dordrecht etc.:
Kluwer Academic Publishers, 584.


\bibitem{Krein}
Krein, M. G., Nudelman, A.A., 1977. The Markov moment problem and extremal problems.
Translations of Mathematical Monographs. Vol. 50. Providence, R.I.: American Mathematical Society.


\bibitem{Luz_Mokl_book2019}
 Luz, M. M., Moklyachuk, M. P., 2019.
 Estimation of stochastic processes with stationary increments and cointegrated sequences.
 London: ISTE; Hoboken, NJ: John Wiley \& Sons, 282 p., 2019.


\bibitem{Luz_Mokl_book2024}
 Luz, M. M., Moklyachuk, M. P., 2024.
Non-stationary stochastic processes estimation: vector stationary increments, periodically stationary multi-seasonal increments.
Berlin: De Gruyter, 292.



\bibitem{Moklyachuk:1994}
Moklyachuk, M. P., 1994.
Stochastic autoregressive sequences and minimax interpolation.
Theory Probab. Math. Stat. 48, 95-103.

\bibitem{Moklyachuk:2000}
Moklyachuk, M. P., 2000.
Robust procedures in time series analysis.
Theory Stoch. Process. 6(3-4), 127-147.

\bibitem{Moklyachuk:2001}
Moklyachuk, M. P., 2001.
Game theory and convex optimization methods in robust estimation problems.
Theory Stoch. Process. 7(1-2), 253-264.

\bibitem{Moklyachuk:2008}
Moklyachuk, M. P., 2008a.
Robust estimations of  functionals of stochastic processes, Kyiv University, Kyiv, 320.

\bibitem{Moklyachuk:2008nonsm}
Moklyachuk, M. P., 2008b.
Nonsmooth analysis and optimization, Kyiv University, Kyiv, 400.

\bibitem{Golichenko}
 Moklyachuk, M.P., Golichenko, I.I.,   2016.
 Periodically Correlated Processes Estimates, LAP
Lambert Academic Publishing, 308.

\bibitem{Moklyachuk:2012}
Moklyachuk, M., Masyutka, O., 2012.
Minimax-robust estimation technique for stationary stochastic processes, LAP
Lambert Academic Publishing, 296.


\bibitem{MoklyachukOstapenko1}
 Moklyachuk M. P., Ostapenko V. I., 2015.
Minimax interpolation problem for harmonizable stable sequences with noise observations.
Journal of Applied Mathematics and Statistics 2(1), 21--42.
doi:10.7726/jams.2015.1003

\bibitem{MoklyachukOstapenko2}
Moklyachuk M. P., Ostapenko V. I., 2016.
Minimax interpolation of harmonizable sequences.
Theor. Probability and Math. Statist. 92, 135-146.

\bibitem{Pshenychn}
Pshenichnyj, B.N., 1971. Necessary conditions for an extremum. Pure and Applied mathematics. 4. New York: Marcel Dekker, 230.


\bibitem{Pourahmadi1984}
Pourahmadi, M., 1984.
On minimality and interpolation of harmonizable stable processes.
SIAM J. Appl. Math. 44(5), 1023-1030.


\bibitem{Rockafellar}
Rockafellar, R. T., 1997.
Convex Analysis. Princeton University Press, 451.


\bibitem{Singer}
Singer, I., 1970.
Best Approximation in Normed Linear Spaces by Elements of Linear Subspaces. Berlin-Heidelberg-New York: Springer-Verlag, 415.

\bibitem{Vastola}
Vastola, K. S., Poor, H. V., 1983. An analysis of the effects of spectral uncertainty on Wiener filtering.
Automatica 28, 289-293.

\bibitem{Weron}
Weron, A., 1985. Harmonizable stable processes on groups: spectral, ergodic and interpolation properties.
Z. Wahrsch. Verw. Gebiete 68(4), 473--491.






\end{thebibliography}
\end{document}